\documentclass[12pt]{amsart}

\usepackage{graphicx}
\usepackage{caption}
\usepackage{subcaption}
\usepackage{color}

\theoremstyle{definition}

\theoremstyle{remark}

\newcommand{\R}{\mathbb{R}}  

\begin{document}


\title[FSI Simulation with Hyperelastic Models]{Numerical Simulation of Fluid-Structure Interaction Problems with Hyperelastic Models I: A Partitioned Approach}

\author{Ulrich Langer}
\address{Johann Radon Institute for Computational and Applied Mathematics (RICAM), Austrian Academy of Sciences and Institute of Computational Mathematics, Johannes Kepler University, Altenberger Strasse 69, A-4040 Linz, Austria}
\email{ulanger@numa.uni-linz.ac.at}
\urladdr{http://www.numa.uni-linz.ac.at/~ulanger/}


\author{Huidong Yang}
\address{Johann Radon Institute for Computational and Applied Mathematics (RICAM), Austrian Academy of Sciences, Altenberger Strasse 69, A-4040 Linz, Austria}
\email{huidong.yang@oeaw.ac.at}
\urladdr{http://people.ricam.oeaw.ac.at/h.yang/}



\begin{abstract}
  In this work, we consider fluid-structure interaction simulation 
  with nonlinear hyperelastic models in the solid part. 
  We use a partitioned 
  approach to deal with the coupled 
  nonlinear fluid-structure interaction problems. We focus on 
  handling the nonlinearity of the fluid and structure sub-problems, 
  the near-incompressibility of materials, the stabilization of 
  employed finite element discretization, 
  and the robustness and efficiency of Krylov subspace 
  and algebraic multigrid 
  methods for the linearized algebraic equations. 
\end{abstract}


\maketitle








\section{Introduction}
Recently, some nonlinear hyperelastic models 
have been adopted in the fluid-structure interaction (FSI) simulation; 
see, e.g., \cite{GJF12:00, YB11:00, Zhang12:00, TTE08:00, YB10:00}. 
In contrast to linear, isotropic and elastic models subjected to 
large displacements and small strains in the FSI simulation 
(see, e.g.,\cite{SC:05, SP11:00}), 
the nonlinear hyperelastic models are often used to 
describe materials undergoing finite deformation and very large 
strains, that are characterized by the anisotropic 
effects, material coefficient jumps across distinct layers, and 
near-incompressibility constraints; 
see, e.g., the arterial layers modeling 
in \cite{Holzapfel00:00, Holzapfel06:00} among others. As phrased 
in \cite{YB10:00}, when incorporating the above characteristics is not 
feasible, some simplified isotropic solid models 
may be used to 
predict wall deformation under the action of the fluid forces. 
In this work, we use a model of modified hyperelastic 
Mooney-Rivlin material (see, e.g., \cite{JB08:00}), 
and a model of thick-walled artery 
with the media and adventitia layer 
(see \cite{Holzapfel00:00}). For this 
two-layer thick-walled artery, we assume 
the angles between collagen fibers and 
the circumferential direction in the media and adventitia are specified for 
a healthy young arterial wall. In \cite{GPlank12:00}, a novel 
Laplace-Dirichlet Rule-Based (LDRB) algorithm to assign fiber 
orientation to heart models has been presented. Besides, 
there are other important issues concerning the arterial modeling, 
e.g., how to obtain arterial wall prestress and to incorporate 
it into the solid momentum balance equation (see \cite{YB11:00}) or 
how to estimate the element-based zero-stress state for arterial 
FSI computations in order to obtain correct reference geometry 
for the finite hyperelastic models (see \cite{TET13:00}), and 
how to derive 
adequate boundary conditions for the solid domain due to 
the external tissue support surrounding 
the arterial wall (see, e.g., \cite{GJF12:00, GJF13:00}). 
These issues are beyond 
the scope of this work, we thus refer 
readers to the corresponding references. 
   
The material coefficient jumps of the two-layer thick-walled artery 
prescribed in \cite{Holzapfel00:00} are not so high, that is 
less problematic for our numerical simulation. Therefore, we concentrate on 
the near-incompressibility of the materials with very large bulk modulus 
$\kappa$. This brings difficulties to 
develop stable discretization and efficient iterative 
solution methods. To enforce incompressibility, 
in \cite{Holzapfel02:00}, a displacement-dilatation-pressure three-field 
Hu-Washizu variational approach combined with 
an augmented Lagrangian optimization technique has been employed. 
In \cite{CA13:00}, a finite element tearing and 
interconnecting (FETI, see, e.g., \cite{PC13:00}) 
method is used to simulate 
biological tissues characterized by anisotropic and nonlinear materials. 
However a pure displacement formulation hired therein may lead to 
{\it locking phenomena}. 

In this work, we use a displacement-pressure 
two-filed variational approach to overcome the {\it locking phenomena}. 
On the one hand, the stabilized equal order mixed 
finite element method in hyperelastic finite deformation 
using a residual-based 
stabilized formulation is well understood, as described in 
\cite{Klaas96:00, Maniatty02:00} for two different types of 
Neo-Hooke material models. An application of this approach to 
the nonlinear inverse problem for hyperelastic tissues has been 
presented in \cite{Goenezen11:00}. 
This mixed displacement-pressure method 
was also used in the FSI simulation of a patient-specific atherosclerotic 
coronary artery modeled as hyperelastic (Fung) material; 
see \cite{TTE09:00}. For the linear 
elastic model, the error estimates of the mixed 
pressure-displacement formulation 
have been analyzed in \cite{SR91:00}. The application to the FSI 
simulation using such a mixed formulation 
has been reported in \cite{YH11:00}. 

On the other hand, the linear system of algebraic 
equations arising from the finite element 
discretization of the mixed 
hyperelastic formulation 
(after linearization using Newton's method \cite{PD05:00}) 
has a very nice form suitable 
for a (nearly) optimal algebraic multigrid (AMG) method. 
This AMG method has been utilized to solve the 
discretized {\it Oseen} and linear elastic equations in 
the FSI simulation; see \cite{YH11:00}, where 
a stabilized coarsening strategy and an 
efficient {\it Braess-Sarazin} smoother 
(see the original and approximate versions 
in \cite{Braess97:00} and \cite{WZ00:00}, respectively) have been 
applied to both fluid and structure sub-problems. 
The number of iterations for the AMG method to converge is 
independent of the problem size and the incompressibility. 
The original contribution of such an AMG method for the fluid problem 
({\it Oseen} equations) 
can be found in, e.g., \cite{WM04:00, WM06:00}. Therein, 
in order to keep $\inf-\sup$ stability on coarse levels 
(since meshes on coarse levels are not available), 
a special coarsening strategy with a proper 
scaling technique is discussed. The two-grid convergence 
of such an AMG method for the Stokes equations 
using a \textit{Petrov-Galerkin} type 
coarse grid operator is shown in \cite{BM:13}. 
 
For the hyperelastic problems, unfortunately a direct application of such 
an AMG method fails. It is found from 
our numerical experiments, that the {\it Braess-Sarazin} smoother does not 
damp efficiently high-frequency errors. We therefore consider 
another class of commonly used 
smoothers based on corrections of many small local problems, 
that are projected from the global problem on properly chosen local patches, 
in a Gauss-Seidel manner. 
For the geometrical multigrid (GMG) methods in 
computational fluid mechanics, the so-called {\it Vanka} smoother based on 
a symmetrical coupled Gauss-Seidel (SCGS) was tested in \cite{Vanka86:00}, 
where a finite-difference formulation using 
staggered locations of the velocity 
and pressure variables is employed. 
An additive version of such a multiplicative 
{\it Vanka}-type smoother under the finite element 
discretization has been analyzed in \cite{ZW03:00}. The {\it Vanka} smoother 
has been widely used in computational fluid dynamics; 
see, e.g., \cite{John01:00, John06:00}, where some stable finite 
element pairs with higher order velocity spaces are used. However, according 
to \cite{WM04:00}, the {\it Vanka} smoother 
deteriorates rapidly in the AMG method for solving 
the 3D {\it Oseen} equations discretized 
with the stabilized $P_1-P_1$ pair. 
This is also observed in our numerical experiments to 
solve the linearized Navier-Stokes (NS) equations in a Newton iteration, 
especially for high Reynolds number.  

The situation is different in the hyperelastic models. We 
apply this smoother to the multigrid method for our hyperelastic problems. 
Using a proper relaxation, the smoother works well under 
sufficient smoothing steps. In \cite{TS09:00}, it 
was observed, the {\it Vanka} smoother in 
the GMG methods works quite robustly 
for the $Q_1-Q_1$ discretized finite 
elastic problem using a hyperelastic Neo-Hooke material in 2D. 
Compared with the {\it Braess-Sarazin} smoother, 
in each {\it Vanka} smoothing step, 
there appear many small local sub-problems, that are usually tackled by 
direct solvers. Thus altogether it leads to a complexity approximately 
proportional to $Nn^3$, where $N$ denotes the number of patches, and $n$ 
the average local sub-problem size. 

We point out, that the 
AMG method we used in this work is not aiming to 
compete with the FETI method in, e.g., \cite{CA13:00}, since parallelization 
has been applied there. However, the (nearly) optimal AMG method 
may be used to solve sub-domain problems in the FETI framework. 
As a comparison, we also consider some Krylov subspace methods 
combined with efficient 
preconditioners, that are derived from the block $LU$ 
factorization of the mixed 
algebraic equations (see \cite{PSV08:00}). 
The iteration numbers of these Krylov subspace methods 
slightly increase when the mesh is refined, that usually 
leads to suboptimal solvers. 

The reminder of this paper is organized in the following way. 
In Section \ref{sec:pre}, we describe the preliminary of this work. 
The complete coupled FSI system is formulated in Section \ref{sec:cs}. 
Section \ref{sec:dis} deals with the temporal and spatial discretization 
of the fluid and structure sub-problems, and their linearization using 
Newton's method. In Section \ref{sec:sm}, solution methods for the nonlinear 
FSI coupled problem, and the nonlinear fluid and structure sub-problems 
are discussed. Some numerical experiments 
are presented in Section \ref{sec:ne}. Finally, the conclusions are 
drawn in Section \ref{sec:con}. 

\section{Preliminary}\label{sec:pre}
\subsection{A model problem}
\begin{figure}[htbp]
  \centering
  \begin{subfigure}[h]{0.4\textwidth}
    \centering
    \includegraphics[scale = 0.3]{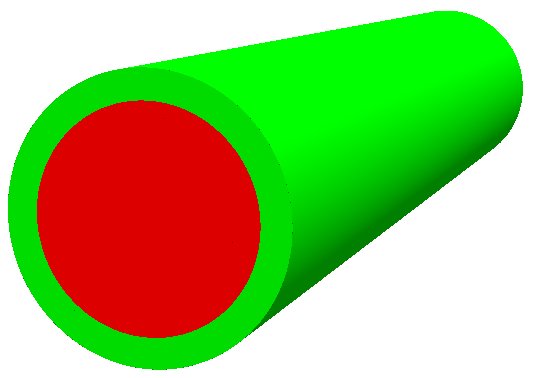}
    \caption{An one-layered model}\label{fsigeoiso}
  \end{subfigure}
  \quad
  \begin{subfigure}[h]{0.4\textwidth}
    \centering
    \includegraphics[scale = 0.3]{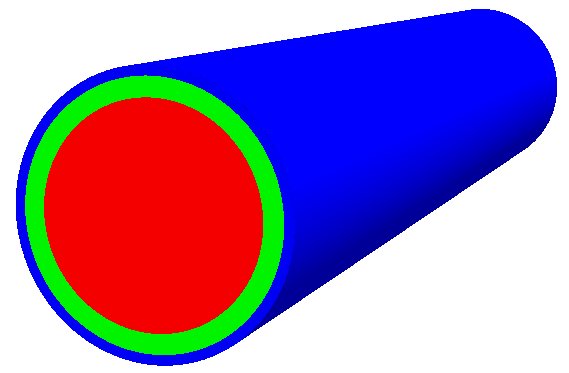}
    \caption{A two-layered model}\label{fsigeoaniso}
  \end{subfigure}
  \begin{subfigure}[h]{0.8\textwidth}
    \centering
    \includegraphics[scale = 1.0]{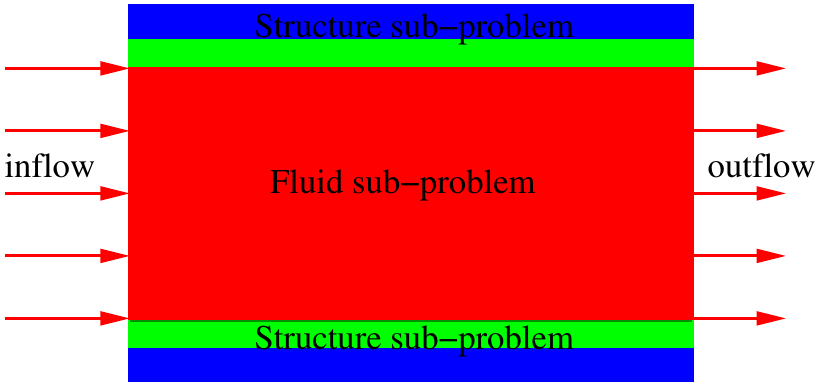}
    \caption{A FSI model problem consisting of fluid and structure sub-problems}\label{fsimodel}
  \end{subfigure}
  \caption{FSI with one-layered (isotropic) and two-layered (anisotropic) hyperelastic models}\label{fsigeomodel}
\end{figure}

As an illustration in Fig. \ref{fsigeomodel}, 
a fluid-structure interaction (FSI) driven by 
an inflow condition is considered, 
which is decomposed into fluid and structure sub-problems, 
and an interaction in between. Considering different application fields, 
we include both isotropic and anisotropic hyperelastic models. 
In the isotropic case, the material holds 
the same properties through the layer; while in the anisotropic case, 
the collagenous matrix of the material varies across the layers.   

\subsection{Geometrical description and arbitrary Lagrangian Euler mapping}
\begin{figure}[htbp]
  \centering
  \scalebox{0.54}{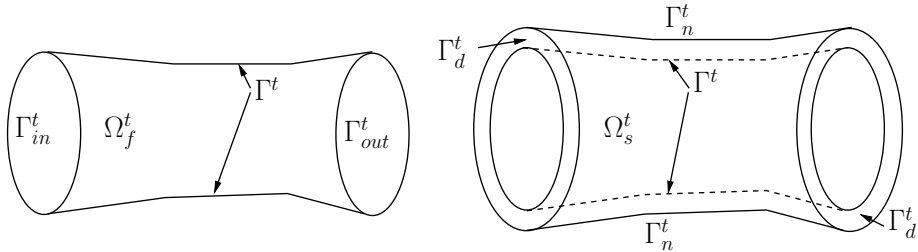}
  \caption{Computational fluid and structure sub-domains obtained using 
    Lagrangian and Arbitrary Lagrangian Eulerian mapping}
  \label{comdm}
\end{figure}

The computational FSI domain $\overline{\Omega^t}\subset\R^3$ 
at time $t$ is composed of 
the fluid sub-domain $\overline{\Omega^t_f}\subset\R^3$ 
and structure sub-domain $\overline{\Omega^t_s}\subset\R^3$: 
$\overline{\Omega^t}=\overline{\Omega^t_s}\cup\overline{\Omega^t_f}$, 
${\Omega^t_s}\cap{\Omega^t_f}=\O$; 
see a schematic representation in Fig. \ref{comdm}. At time $t=0$, 
$\overline{\Omega^0}$ is referred to as a reference (initial) configuration: 
$\overline{\Omega^0}=\overline{\Omega^0_s}\cup\overline{\Omega^0_f}$, 
$\Omega^0_s\cap\Omega^0_f=\O$. 
For the fluid sub-problem, 
$\Gamma_{in}^t$ and $\Gamma_{out}^t$ are parts of boundary $\partial\Omega^t_f$ 
where certain Neumann forces are applied. In an analogous manner, 
$\Gamma_{d}^0$ and $\Gamma_{n}^0$ are employed to denote 
parts of boundary $\partial\Omega^0_s$, where structure Dirichlet and Neumann 
conditions are applied, respectively. The intersection $\Gamma^t$ 
between two sub-domains at time $t$ is called the FSI interface: 
$\Gamma^t=\partial\Omega^t_f\cap\partial\Omega^t_s$. 
Again at $t=0$, $\Gamma^0$ is the reference (initial) interface: 
$\Gamma^0=\partial\Omega^0_f\cap\partial\Omega^0_s$, where certain 
interface conditions are fulfilled. 

In order to track particle motion of the structure body $\Omega_s^0$, 
we adopt the Lagrangian mapping $\mathcal{L}^t (\cdot) : 
\mathcal{L}^t(x_0)=x_0+\hat{d}_s(x_0, t)$ for all $x_0\in\Omega_s^0 $ and 
$t\in (0, T)$, with the structure displacement 
$\hat{d}_s(\cdot, \cdot) : \Omega_s^0\times (0, T]\mapsto\R^3$. 
An arbitrary Lagrangian Euler (ALE, see, e.g., \cite{LF99:00, JD04:00, TH:81}) 
mapping $\mathcal{A}^t(\cdot)$ 
is used to capture the movement of the fluid sub-domain $\Omega_f^0$: 
$\mathcal{A}^t(x_0)=x_0+\hat{d}_f(x_0, t)$ for all 
$x_0\in\Omega_f^0 $ and $t\in (0, T)$, where the fluid displacement 
$\hat{d}_f(\cdot, \cdot): \Omega_f^0\times (0, T]\mapsto\R^3$, only follows 
particle motion of the fluid body on the interface $\Gamma^0$, but is 
arbitrarily extended into the domain $\Omega_f^0$. 
A classical option for the ALE mapping is the so-called harmonic extension: 
\begin{subequations}\label{eq:harext}
  \begin{eqnarray}
    -\Delta \hat{d}_f=0&{\textup{ in }} \Omega_f^0,\\
    \hat{d}_f=\hat{d}_s&{\textup{ on }} \Gamma^0
  \end{eqnarray}
\end{subequations}
with homogeneous Dirichlet boundary conditions 
$\hat{d}_f=0$ on $\Gamma_{in}\cup\Gamma_{out}$. 
The fluid sub-domain velocity 
$w_f:\Omega_f^t\mapsto\R^3$ at time $t$ is then given by
\begin{equation}\label{fdv}
  w_f=\frac{\partial \hat{d}_f}{\partial t}\circ {\mathcal{A}^t}^{-1}.
\end{equation}
In our situation, the simple harmonic extension (\ref{eq:harext}) 
provides satisfying results. 
Other more advanced mesh motion techniques under 
the ALE framework can be found in, e.g., \cite{Wick11:00}.

\subsection{Nonlinear hyperelastic modeling}
In order to formulate the hyperelastic
sub-problem in Section \ref{sec:strucsub}, we 
first introduce the following useful notations in nonlinear 
continuum mechanics, that can be found in, 
e.g., \cite{GAH00:00, JB08:00, RW97:00, WP:06}. 
Let $F$ denote the deformation gradient tensor 
given by $F=\partial\mathcal{L}^t/\partial x_0=I+\nabla\hat{d}_s$ and $J$ 
the determinant given by $J=\textup{det}F$. We then define the right 
Cauchy-Green tensor $C$ as $C=F^TF$ and the three 
principal invariants are respectively given by 
\begin{equation}
I_1=C:I, \quad I_2=\frac{1}{2}(I_1^2-C:C), \quad I_3=\textup{det}(C).
\end{equation}
For a hyperelastic material, the second Piola-Kirchoff tensor $S$ 
is defined as 
\begin{equation}
S=2\frac{\partial \Psi}{\partial C},
\end{equation}
where $\Psi$ represents the energy functional 
depending on the invariants. The energy functional 
is chosen according to material properties, and the 
second Piola-Kirchoff tensor is then computed accordingly.

\subsubsection{A modified model of Mooney-Rivlin material}
We first consider a modified model of 
Mooney-Rivlin material, for which the energy functional is given by
\begin{equation}\label{eq:mmreng}
\Psi=\frac{c_{10}}{2}(J_1-3)+\frac{c_{01}}{2}(J_2-3)+\frac{1}{2}\kappa(J-1)^2,
\end{equation}
where the first two invariants are respectively given by 
\begin{equation}
J_1=I_1I_3^{-1/3},\quad J_2=I_2I_3^{-2/3}. 
\end{equation}
Material parameters $c_{10}>0$ and $c_{01}>0$ are related to shear 
modulus $\mu^l$ for consistency with linear elasticity in the 
limit of small strains: $\mu^l=2(c_{10}+c_{01})$. The bulk modulus $\kappa>0$ 
indicates material compressibility: the larger the bulk modulus is, 
the more incompressible the material becomes. 
In this case, the second Piola-Kirchoff tensor is given by 
\begin{equation}\label{eq:mmrpkt}
  S=c_{10}\frac{\partial J_1}{\partial C}
  +c_{01}\frac{\partial J_2}{\partial C}+2
  \kappa(J-1)\frac{\partial J}{\partial C}.
\end{equation} 
In order to overcome the usual {\it locking phenomena} 
with large bulk modulus, we 
introduce the pressure $\hat{p}_s:=\hat{p}_s(x, t)
=-\kappa(J-1):\Omega_s^0\times (0, T]\mapsto \R^3$ for all $x\in\Omega_s^0$ 
at time $t$. Using the notation 
$$S^{'}=c_{10}\frac{\partial J_1}{\partial C}
+c_{01}\frac{\partial J_2}{\partial C},$$ 
then the second Piola-Kirchoff tensor is rewritten as 
\begin{equation}\label{eq:mmrpktmix}
S=S^{'}-\hat{p}_sJC^{-1},
\end{equation}
where we use $\partial J/\partial C = JC^{-1}/2$. 
The Mooney-Rivlin model is often used 
to model elastic response of rubber-like 
materials (see, e.g., \cite{JB08:00}). 
We use this model in the FSI 
simulation as an isotropic 
hyperelastic case for test purpose; see the one-layered model in 
Fig. \ref{fsigeomodel}. 
\subsubsection{A two-layer thick-walled artery}
In many applications, simple isotropic 
hyperelastic models might be insufficient to describe 
the mechanical response of materials 
like arterial tissues; 
see, e.g., \cite{Holzapfel00:00, Holzapfel06:00, CA13:00}. 
When studying mechanical response of a healthy young arterial (with no 
pathological intimal changes) in hemodynamics, we usually 
model the artery as a two-layer thick-walled tube consisting of 
the media and the adventitia layers. 
The media layer has a common FSI 
interface with the fluid. The energy functional of such an arterial model 
is prescribed by:
\begin{equation}\label{eq:areng}
  \Psi=\Psi_{\text{iso}}+\Psi_{\text{aniso}}+\frac{1}{2}\kappa(J-1)^2,
\end{equation}    
where the isotropic response in each layer and the strain 
energy stored in the collagen fibers are respectively prescribed 
by the classical neo-Hooken model and an exponentional function 
(see \cite{Holzapfel00:00,GAH06:00}): 
\begin{equation}
  \Psi_{\text{iso}}=\frac{c_{10}}{2}(J_1-3),\quad 
  \Psi_{\text{aniso}}=\frac{k_1}{2k_2}\displaystyle\sum_{i=4,6}(\exp(k_2(J_i-1)^2)-1).
\end{equation}    
Here $c_{10}>0$ is a stress-like material parameter 
and $\kappa>0$ the bulk modulus as in the 
modified Mooney-Rivlin material, $k_1>0$ a stress-like material parameter, 
and $k_2>0$ a dimensionless parameter. The invariants $J_4>1$ 
and $J_6>1$ (active in extension) are given 
by:
\begin{equation}\label{eq:arinv46}
    J_4=I_3^{-1/3}A_1:C,\quad J_6=I_3^{-1/3}A_2:C,
\end{equation}
where the tensors $A_1$ and $A_2$ characterizing 
the media and adventitia structure are prescribed by 
\begin{equation}
A_1=a_{01}\otimes a_{01},\quad A_2=a_{02}\otimes a_{02}
\end{equation}
with the direction vectors $a_{01}$ and $a_{02}$ that are specified by 
the angle $\alpha$ between the collagen fibers and 
the circumferential direction in the media or adventitia. 
In particular,  they have the forms 
\begin{equation}
a_{01}=(0, \cos\alpha, \sin\alpha)^T,\quad 
a_{01}=(0, \cos\alpha, -\sin\alpha)^T. 
\end{equation}
Then, the second Piola-Kirchoff tensor $S$ for this hyperelastic model 
is computed as follows:
\begin{equation}\label{eq:arpkt}
  S=c_{10}\frac{\partial J_1}{\partial C}+
  k_1\displaystyle\sum_{i=4,6}(\exp(k_2(J_i-1)^2)(J_i-1)
  \frac{\partial J_i}{\partial C})+
  2\kappa(J-1)\frac{\partial J}{\partial C}.
\end{equation} 
Using the pressure $\hat{p}_s$ as defined for 
the modified Mooney-Rivlin material and the notation 
$$S^{'}=c_{10}\frac{\partial J_1}{\partial C}+
k_1\displaystyle\sum_{i=4,6}(\exp(k_2(J_i-1)^2)(J_i-1)
\frac{\partial J_i}{\partial C}),$$ 
the second Piola-Kirchoff tensor is reformulated 
in the following form
\begin{equation}\label{eq:arpktmix}
   S=S^{'}-\hat{p}_sJC^{-1}.
\end{equation}
The material parameters $c_{10}$, $k_1$, $k_2$, and the 
angle $\alpha$ in the model (\ref{eq:areng}) take 
different values in the media and adventitia, i.e., they 
have jumps across the layers. In particular, the constants $k_1$ and 
$k_2$ are associated with the anisotropic contribution of collagen 
to the response; see an illustration in 
Fig. 14 of \cite{Holzapfel00:00}, from which 
we also adopt the material and geometrical data of 
a carotid artery from a rabbit as a benchmark in our numerical experiments. 
Therefore, we consider 
this model as an anisotropic hyperelastic model in the FSI simulation; 
see the illustration of two-layered model in Fig. \ref{fsigeomodel}.

\subsubsection{A short remark on hyperelastic modeling}
As stated in \cite{Holzapfel00:00}, although the material parameters 
are constants independent of the geometry, 
opening angle or fiber angle, the internal pressure/radius response does 
depend on these quantities. For simplicity, 
we use the same assumption as therein 
that the fiber angles are the same in the load-free configuration. 
Unfortunately, in patient-specific arterial FSI computations, 
the image-based arterial geometry is usually not stress-free. 
In order to overcome this difficulty, in, e.g., \cite{YB11:00}, the 
prestress is directly integrated into the momentum equation, 
or in, e.g., \cite{TET13:00}, the authors 
developed a methodology to find a zero-stress configuration for the 
artery, that is used as a reference configuration for the 
hyperelastic models. Another interesting issue is related to 
proper descriptions of boundary conditions for the hyperelastic model since 
the artery is surrounded by tissues. 
In \cite{GJF12:00, GJF13:00, CNM12:00}, the 
viscoelastic support conditions along the artery wall are proposed 
to model the effect of external tissues. These important issues are not 
the focus of this work. However, the solution methods we are proposing 
are not only restrict to these two hyperelastic models used in this work. 
Indeed it is easy to integrate the above effects into our methodology. 

\section{The coupled system}\label{sec:cs}
\subsection{The fluid sub-problem}
The fluid sub-problem in a deformable domain $\Omega_f^t\subset\R^3$ 
is governed by the incompressible NS equations 
under the ALE framework: Find the fluid velocity $u:=u(x, t):\Omega_f^t\mapsto \R^3$ 
and pressure $p_f:=p_f(x, t):\Omega_f^t\mapsto \R$ 
for all $x\in \Omega_f^t$ at time $t$ such that 
\begin{subequations}\label{eq:ns}
  \begin{eqnarray}
  \rho_f\partial_t u|_{{\mathcal A}^t}+\rho_f((u- w_f)\cdot\nabla) u-\nabla\cdot \sigma_f(u, p_f)=0&{\textup{ in }}\Omega_f^t,\\
  \nabla\cdot u=0&{\textup{ in }}\Omega_f^t
  \end{eqnarray}
\end{subequations}
with Neumann boundary conditions 
\begin{subequations}\label{eq:nsbnd}
  \begin{eqnarray}
  { \sigma}_{f}(u, p_f)n_f=g_{in}&{\textup{ on }}\Gamma_{in}^t,\\
  { \sigma}_{f}(u, p_f)n_f=0&{\textup{ on }}\Gamma_{out}^t,
  \end{eqnarray}
\end{subequations}
where $\rho_f$ denotes the fluid density, $n_f$ the 
outerward unit normal vector 
to $\Gamma_{in/out}^t$, 
$\sigma_{f}(u, p_f):=2\mu\varepsilon(u)-p_fI$ the Cauchy stress tensor, 
$\varepsilon(u):=\frac{1}{2}({\nabla u+(\nabla u)^{T}})$ 
the strain rate tensor, 
$\mu$ the dynamic viscosity, $g_{in}$ the 
given data for the inflow boundary condition, and 
\begin{equation}
\partial_t u|_{{\mathcal A}^t}:=\partial_t u+(w_f\cdot\nabla) u 
\end{equation} 
the ALE time derivative with the fluid sub-domain velocity 
$w_f$ defined by (\ref{fdv}). 

\subsection{The structure sub-problem}\label{sec:strucsub}
For finite deformation, 
the structure sub-problem under the Lagrangian framework 
in a reference domain $\Omega_s^0\subset {\mathbb R}^3$ 
is governed by the following equations: 
Find the structure displacement $\hat{d}_s$ and pressure $\hat{p}_s$ 
for all $x\in \Omega_s^0$ at time $t$ such that 
\begin{subequations}\label{eq:elas}
  \begin{eqnarray}
    \rho_s\partial_{tt}\hat{d}_s- 
    \nabla\cdot (FS)=0&{\textup{ in }}\Omega_s^0,\\  
    \hat{p}_s=-\kappa(J-1)&{\textup{ in }}\Omega_s^0
  \end{eqnarray}
\end{subequations}
with Dirichlet and Neumann boundary conditions
\begin{subequations}\label{eq:elasbnd}
  \begin{eqnarray} 
     \hat{d}_s=0&{\textup{ on }}\Gamma_d^0,\\
     FS\hat{n}_s=0&{\textup{ on }} \Gamma_n^0,
   \end{eqnarray}
\end{subequations}
where $\rho_s$ denotes the structure density, $\hat{n}_s$ the outerward unit normal vector 
to $\Gamma_n^0$ of the initial configuration, and 
the second Piola-Kirchhoff tensor $S$ is given by 
(\ref{eq:mmrpktmix}) or (\ref{eq:arpktmix}) according to the choice of 
the hyperelastic models. 
\subsection{Interface equations}
The interface equations on $\Gamma^0$ 
are respectively described by the 
classical no-slip condition and equivalence of 
surface tractions:
\begin{subequations}\label{eq:interface}
  \begin{eqnarray}
    u\circ\mathcal{A}^t=\partial_t \hat{d}_s&{\textup{ on }} \Gamma^0,\\
    \sigma_{f}(u, p_f)n_f\circ \mathcal{A}^t+FS\hat{n}_s=0&{\textup{ on }}\Gamma^0.
  \end{eqnarray}
\end{subequations}
The fully coupled FSI problem consists of (\ref{eq:harext}), 
(\ref{eq:ns})-(\ref{eq:interface}), and proper initial conditions: 
$u(x, 0)=0$ for all $x\in\Omega_f^0$, $\hat{d}_s(x, 0)=\partial_t\hat{d}_s(x, 0)=0$ 
for all $x\in\Omega_s^0$.

\section{Temporal and spatial discretization and linearization}\label{sec:dis}
\subsection{Temporal discretization}
In order to solve the FSI problem, we 
first make temporal discretization and seek 
the FSI solution at each time level. Thus we subdivide 
the time period $(0, T]$ into $N$ equidistant intervals and 
the time step size is given by $\Delta t=T/N$. 
We aim to find the FSI solution at each time 
level $t^n=n\Delta t$, $n=1, ..., N$. 
At time level $t^0$, the FSI solution is known by the initial conditions. 
For simplicity of notations, a function $f$ at time level $t^n$ is denoted 
by $f^n:=f(\cdot, t^n)$.
\subsubsection{Temporal discretization for the fluid sub-problem}
An fully implicit Euler scheme for the fluid sub-problem is adopted: 
\begin{equation}\label{eq:nstdis}
\partial_t u (x, t^n)|_{{\mathcal A}^{t^n}}\approx (u^n-u^{n-1}\circ{\mathcal A}^{t^{n-1}}\circ{({\mathcal A}^{t^n})}^{-1})/\Delta t.
\end{equation}
At each nonlinear FSI iteration, 
the fluid sub-domain $\Omega_f^{t^n}$ at time level $t^n$ is extrapolated 
by the sub-domain from the previous iteration. The fluid domain velocity $w_f^t$ 
in the convection term is then computed using the first order Euler scheme:
\begin{equation}\label{eq:hartdis}
  w_f^n=(\tilde{d}_f^n-\hat{d}_f^{n-1})\circ{({\mathcal A}^{t^n})}^{-1}/\Delta t,
\end{equation}  
where $\tilde{d}_f^n$ is computed by solving the harmonic 
extension problem (\ref{eq:harext}) with obtained interface displacement from 
previous nonlinear iteration, and $\hat{d}_f^{n-1}$ 
is the fluid artificial domain displacement from previous time level $t^{n-1}$. 

By the ALE property, the fluid velocity $u^n_{|\Gamma^t}$ 
on the interface is then given by the fluid 
domain velocity, i.e., $u^n_{|\Gamma^t}=w^n_{|\Gamma^t}$. In addition, since 
we use a fully implicit scheme, the nonlinearity in the 
fluid convection term $u^n\cdot\nabla u^n$ needs to be 
tackled by Newton's method. 
\subsubsection{Temporal discretization for the structure sub-problem} 
To discretize the structure sub-problem 
in time, a first order Newmark-$\beta$ scheme 
is used; see \cite{NM59:00}. For simplicity of notations, let 
$\dot{d}_s^n:=\partial_{t}\hat{d}_s(x, t^n)$ 
and $\ddot{d}_s^n:=\partial_{tt}\hat{d}_s(x, t^n)$. Then the 
Newmark-$\beta$ scheme is given by
\begin{subequations}
  \begin{eqnarray}
    \ddot{d}_s^{n}&\approx&\frac{1}{\beta\Delta t^2}
    (\hat{d}_s^{n}-\hat{d}_s^{n-1})-
    \frac{1}{\beta\Delta t}\dot{d}_s^{n-1}-
    \frac{1-2\beta}{2\beta}\ddot{d}_s^{n-1},\\
    \dot{d}_s^{n}&\approx&\dot{d}_s^{n-1}
    +\gamma\Delta t\ddot{d}_s^{n}+
    (1-\gamma)\Delta t\ddot{d}_s^{n-1},
  \end{eqnarray}
\end{subequations}
with constants $0<\beta\leq 1$ and $0\leq\gamma\leq 1$. 
As observed in \cite{HY12:00} for 
linear elastic models using 
the mixed displacement and pressure formulation, 
the conventional choice of 
$2\beta=\gamma=1$ in the 
numerical experiments leads to temporal oscillation of the pressure field, and 
the choice of $\beta>0.5$ leads to stable time discretization. 
Thus in our numerical experiments, we have chosen $\beta=0.625$ and $\gamma=1$ 
in order to avoid the pressure oscillation in time. 

\subsection{Time semi-discretized weak formulations}
\subsubsection{Function spaces}
We aim to find weak FSI solution on proper function spaces. 
Let $H^1(\Omega_f^0)$, $H^1(\Omega_s^0)$, $L^2(\Omega_f^0)$ and $L^2(\Omega_s^0)$ 
denote the standard Sobolev and Lebesgue spaces (see, e.g., \cite{AF03:00}) 
on $\Omega_f^0$ and $\Omega_s^0$, respectively. 
The function spaces for the 
fluid velocity and pressure are 
respectively given by 
$V_f^t:=\{v_f : v_f\circ{\mathcal A}^t\in H^1(\Omega_f^0)^3\}$ and 
$Q_f^t:=\{q_f : q_f\circ{\mathcal A}^t\in L^2(\Omega_f^0)\}$; 
the function spaces $V_s$ and $Q_s$ 
for the structure displacement and 
pressure shall be appropriately chosen according to the 
nonlinearities; see, e.g., \cite{JMB:76,PGC88:00}. 

\subsubsection{The time semi-discretized fluid weak formulation}
After time discretization, the weak formulation of the fluid sub-problem reads: Find 
$u^n\in V_{f, D}^{t^n}:=\{v\in V_f^{t^n} : v=g_d\in\R^3 \text{ on }\Gamma^{t^n}\}$ 
and $p_f^n\in Q_f^{t^n}$ such that 
for all $v\in  V_{f, 0}^{t^n}:=\{v\in V_f^{t^n} : v=0 \text{ on }\Gamma^{t^n}\}$ 
and $q\in Q_{f}^{t^n}$
\begin{subequations}\label{eq:weakns}
  \begin{eqnarray}
    \tilde{R}_f^1((u^n, p_f^n), v)=0,&\\
    \tilde{R}_f^2((u^n, p_f^n), q)=0,&
  \end{eqnarray}
\end{subequations}
where the residuals $\tilde{R}_f^1(\cdot, \cdot)$ and $\tilde{R}_f^2(\cdot, \cdot)$ 
on the function spaces $V_f^{t^n}$ and $Q_f^{t^n}$ are given by 
\begin{equation*}
  \begin{aligned}
    \tilde{R}_f^1((u, p), v)&=(\frac{\rho_f}{\Delta t}u
    +\rho_f((u-w^n)\cdot\nabla)u, v)_{\Omega_f^{t^n}}
    +2\mu(\varepsilon(u), \varepsilon(v))_{\Omega_f^{t^n}}\\
    &-(p, \nabla\cdot v)_{\Omega_f^{t^n}}
    -(\frac{\rho_f}{\Delta t}u^{n-1}, v)_{\Omega_f^{t^n}}
    -\langle g_{in}, v\rangle_{\Gamma^{t^n}},\\
    \tilde{R}_f^2((u, p), q)&=-(q, \nabla\cdot u)_{\Omega_f^{t^n}}.
  \end{aligned}
\end{equation*}
Here we assume a Dirichlet boundary condition $g_d$ 
is prescribed on the interface $\Gamma^{t^n}$. 
\subsubsection{The time semi-discretized structure weak formulation}
In an analogous way, the time semi-discretized weak 
formulation of the structure sub-problem reads: 
Find $d_s^n\in V_{s,0}:=\{v\in V_s : v=0 \text{ on }\Gamma_d^0\}$ and $p_s^n\in Q_s$ 
for all $v\in  V_{s, 0}$ and $q\in Q_{s}$ such that 
\begin{subequations}\label{eq:weakelas}
  \begin{eqnarray}
    \tilde{R}_s^1((d_s^n, p_s^n), v)=0,&\\
    \tilde{R}_s^2((d_s^n, p_s^n), q)=0,& 
  \end{eqnarray}
\end{subequations}
where the residuals $\tilde{R}_s^1(\cdot, \cdot)$ and $\tilde{R}_s^2(\cdot, \cdot)$ 
on the function spaces $V_s$ and $Q_s$ are given by
\begin{equation*}
  \begin{aligned}
    \tilde{R}_s^1((d, p), v)&=(\frac{\rho_s}{\beta\Delta t^2}d, v)_{\Omega_s^0}
    +(S^{'}, (F^T\nabla v))_{\Omega_s^0}
    -(pJF^{-T}, \nabla v)_{\Omega_s^0}\\
    &-\langle g_{n}, v\rangle_{\Gamma^0}-(r_s, v)_{\Omega_s^0},\\
    \tilde{R}_s^2((d, p), q)&=-(J-1, q)_{\Omega_s^0}-\frac{1}{\kappa}(p, q)_{\Omega_s^0},
  \end{aligned}
\end{equation*}
where $r_s=\frac{\rho_s}{\beta\Delta t^2}d_s^{n-1}
+\frac{\rho_s}{\beta \Delta t}\dot{d}_s^{n-1}
+\frac{\rho_s(1-2\beta)}{2\beta}\ddot{d}_s^{n-1}$. 
Here we assume a Neumann data $g_{n}\in\R^3$ 
on the interface $\Gamma_0$ is known. 
\subsection{Stabilized finite element methods}
Finite element discretization requires a triangulation of the 
computational FSI domain. For this purpose, we 
use Netgen \cite{JS97:00} to 
generate a tetrahedral mesh with conforming grids 
on the FSI interface. In addition it 
resolves different structure layers with 
prescribed material parameters; 
see an illustration in Fig. \ref{fsimesh}. 
\begin{figure}[htbp]
  \begin{subfigure}[h]{0.45\textwidth}
    \centering
    \includegraphics[scale = 0.25]{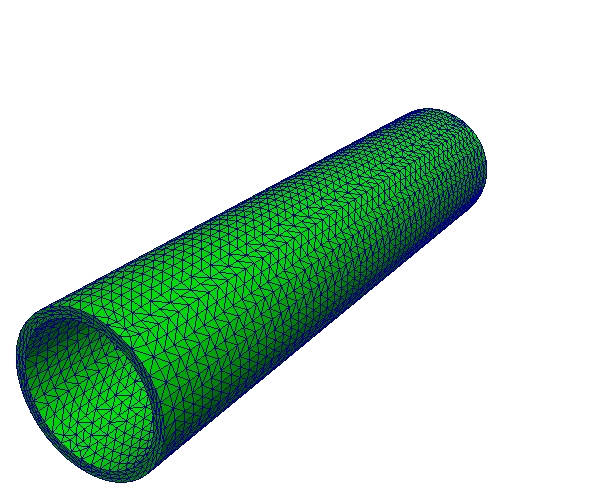}
    \caption{A triangulation of an inner layer of $\Omega_s^0$}\label{insmesh}
  \end{subfigure}
  \begin{subfigure}[h]{0.45\textwidth}
    \centering
    \includegraphics[scale = 0.25]{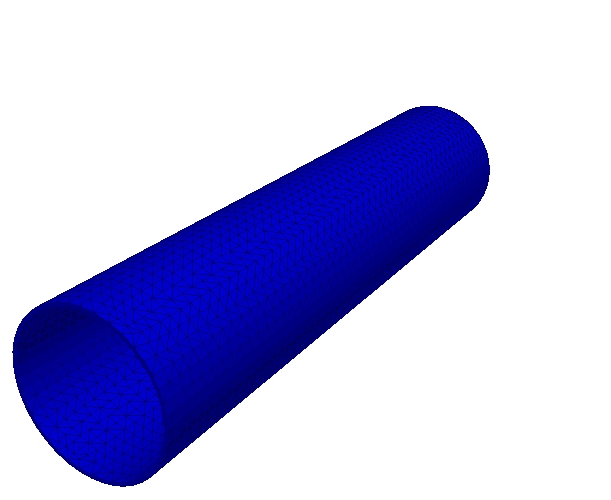}
    \caption{A triangulation of an outer layer of $\Omega_s^0$ }\label{ousmesh}
  \end{subfigure}
  \begin{subfigure}[h]{0.6\textwidth}
    \centering
    \includegraphics[scale = 0.25]{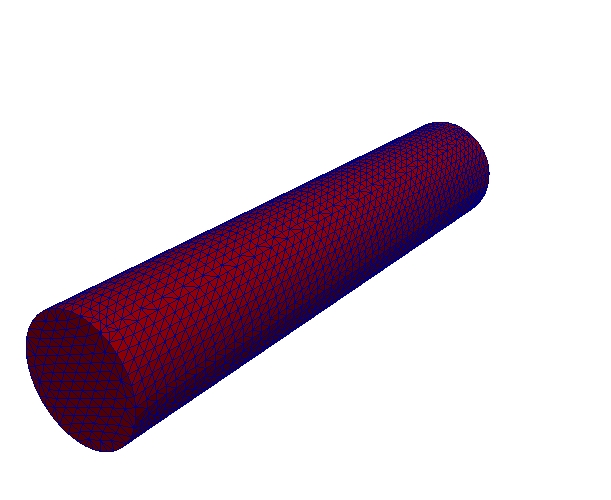}
    \caption{A triangulation of $\Omega_f^0$}\label{fmesh}
  \end{subfigure}
  \caption{An example of a FSI tetrahedral mesh with matching 
    grids on the interface using netgen. It resolves 
    different layers (cf. Fig. \ref{fsigeomodel})}
\end{figure}\label{fsimesh} 

Equal order linear finite element spaces for 
the velocity and pressure fields in the fluid sub-problem, 
and the displacement and pressure fields in the structure 
sub-problem suffer from instability. 
Therefore some stabilization techniques are employed and adapted. 
\subsubsection{Stabilization for the fluid sub-problem}
Let $\mathcal{T}_{f,h}$ be the triangulation of the fluid sub-domain 
$\Omega_f^{t^n}$, that 
is obtained by the ALE mapping on the triangulation  $\mathcal{T}_{f,h}^0$ 
of the reference fluid sub-domain $\Omega_f^0$. The 
finite element spaces of the fluid velocity and pressure are given by 
$V_{f,h}:=\{v_f\in C^0(\overline{\Omega_f^t})^3 : v_f|_T=P_1(T)^3,\; 
\forall T\in \mathcal{T}_{f,h}\}$ and 
$Q_{f,h}:=\{q_f\in C^0(\overline{\Omega_f^0}) : q_f|_T=P_1(T),\; 
\forall T\in \mathcal{T}_{f,h}\}$, respectively, where 
$P_1(T)$ denotes the linear polynomial defined on a tetrahedral element. 

As is well known, the discrete solution may contain unphysical 
oscillations, e.g., spurious pressure modes, 
when equal order finite element spaces are used in the fluid problem, that 
do not fulfill the $\inf-\sup$ or LBB (Ladyshenskaya-Babu\v{s}ka-Brezzi) 
stability condition; see, e.g., \cite{IB73:00, FB74:00, BF91:00}. 
In addition, the instability may occur in advection 
dominated regions of the domain. 
To suppress instability, we employ a unified 
stabilization technique, \textit{streamline-upwind/-} and 
\textit{pressure-stabilizing/Petrov-Galerkin}  (SUPG/PSPG) method, 
developed in 
\cite{Brooks82:00, Hughes89:00, TH86:00, TET92:00, JJD94:00, TET00:00}, 
that is able to 
enhance velocity field stability and meanwhile allows the 
use of equal order finite element spaces for 
both the velocity and pressure fields. 
This stabilization technique was adapted to the fluid problem 
in the FSI simulation; see, e.g., \cite{Wall99, WD06:00, YB06:00}. 
In \cite{Yang11:00}, it was used for the 
stabilized finite element discretization of the \textit{Oseen}-type 
equations on hybrid meshes. Applications of 
more recently developed 
stabilization techniques to the fluid problems in the FSI 
simulation can be found in, e.g., \cite{FLD09:00, CF07:00}.

The fully discretized and stabilized fluid sub-problem 
reads: Find $u_h\in V_{f, h}\cap V_{f,d}^{t^n}$ and $p_{f,h}\in Q_{f,h}$ such that for all 
$v_h\in V_{f,h}\cap V_{f,0}^{t^n}=:V_{f,h}^0$ and $q_h\in Q_{f,h}$
\begin{subequations}\label{eq:stabweakfluid}
  \begin{eqnarray}
    R_f^1((u_h^n, p_{f,h}^n), v_h)=0,&\\
    R_f^2((u_h^n, p_{f,h}^n), q_h)=0,& 
  \end{eqnarray}
\end{subequations}
where $R_f^1(\cdot, \cdot)=\tilde{R}_f^1(\cdot, \cdot)+W_f^1(\cdot, \cdot)$ and 
$R_f^2(\cdot, \cdot)=\tilde{R}_f^2(\cdot, \cdot)+W_f^2(\cdot, \cdot)$. Here we 
include $W_f^1(\cdot, \cdot)$ and $W_f^2(\cdot, \cdot)$ for corresponding 
stabilization terms coming from the SUPG/PSPG method. 

\subsubsection{Stabilization for the structure sub-problem}
Let $\mathcal{T}_{s,h}$ be the triangulation of the 
structure sub-domain $\Omega_s^0$; 
let $V_{s,h}:=\{v_s\in C^0(\overline{\Omega_s^0})^3 : v_s|_T
=P_1(T)^3,\; \forall T\in \mathcal{T}_{s,h}\}\cap V_{s,0}$ and  
$Q_{s,h}:=\{q_s\in C^0(\overline{\Omega_s^0}) : q_s|_T=P_1(T),\; 
\forall T\in \mathcal{T}_{s, h}\}$ represent the finite element spaces of 
the structure displacement and pressure, respectively. 

Thanks to the absence of convection term in the hyperelastic model, 
we only face instability due to equal order interpolation spaces 
for both displacement and pressure fields. 
The PSPG method introduced in \cite{TH86:00} 
is used to circumvent the $\inf-\sup$ 
stability condition; see the application 
in, e.g., \cite{Klaas96:00} to two different 
types of Neo-Hooke material models, and 
\cite{Goenezen11:00} to a soft tissue model. We 
adapt this technique to hyperlastic models used in this work. 
The stabilized formulation is obtained by 
augmenting the Galerkin discretization 
with the following introduced least-square 
finite element stabilization term: 
\begin{equation}\label{eq:elasstterm}
  W_s((d_{s,h}, p_{s,h}), q_h)
  =\displaystyle\sum_{T\in\mathcal{T}_{s,h}}\tau(r_{s,h}
  -\frac{\rho_s}{\beta\Delta t^2}d_{s,h}
  +\nabla\cdot (FS),  F^{-T}\nabla q_h)_{T}
\end{equation}
for $d_{s,h}\in V_{s,h}$ and $p_{s,h}, q_h\in Q_{s,h}$, where 
$r_{s,h}=\frac{\rho_s}{\beta\Delta t^2}d_{s,h}^{n-1}
+\frac{\rho_s}{\beta \Delta t}\dot{d}_{s,h}^{n-1}
+\frac{\rho_s(1-2\beta)}{2\beta}\ddot{d}_{s,h}^{n-1}$ 
corresponds to the discretization of the right hand side of 
the structure momentum equation. The stabilized 
parameter $\tau$ is given by $\tau=\frac{h_T^2}{4\mu^l}$, 
where $h_T$ is the characteristic element length for the 
tetrahedral element $T$. The 
parameter $\mu^l$ is given by 
$\mu^l=2(c_{10}+c_{01})$ and 
$\mu^l=2(c_{10}+\frac{k_1}{2k_2})$ 
for the modified Mooney-Rivlin material 
and two-layer thick-walled artery, respectively. 
Since we are working with linear finite elements, the derivative 
of $S^{'}$ cancel out. Then the above 
stabilization term (\ref{eq:elasstterm}) is simplified as follows:
\begin{equation}\label{eq:stelas}
  \begin{aligned}    
    W_s((d_{s,h}, p_{s,h}), q_h)&=\displaystyle\sum_{T\in\mathcal{T}_{s,h}}
    \tau(r_{s,h}-\frac{\rho_s}{\beta\Delta t^2}d_{s,h}
    -\nabla\cdot (p_{s,h}JF^{-T}),  F^{-T}\nabla q_h)_{T}\\
    &= \displaystyle\sum_{T\in\mathcal{T}_{s,h}}
    \tau(r_{s,h}-\frac{\rho_s}{\beta\Delta t^2}d_{s,h}-JF^{-T}\nabla p_{s,h}, 
    F^{-T}\nabla q_h)_{T}.
  \end{aligned}
\end{equation}
The fully discretized and stabilized structure sub-problem 
reads: Find $d_{s,h}\in V_{s,h}$ and $p_{s,h}\in Q_{s,h}$ such that for all 
$v_h\in V_{s,h}$ and $q_h\in Q_{s,h}$
\begin{subequations}\label{eq:stabweakelas}
  \begin{eqnarray}
    R_s^1((d_{s,h}^n, p_{s,h}^n), v_h)=0,&\\
    R_s^2((d_{s,h}^n, p_{s,h}^n), q_h)=0,& 
  \end{eqnarray}
\end{subequations}
where $R_s^1(\cdot, \cdot)=\tilde{R}_s^1(\cdot, \cdot)$ and 
$R_s^2(\cdot, \cdot)=\tilde{R}_s^2(\cdot, \cdot)+W_s(\cdot, \cdot)$. Here 
$W_s$ corresponds to an additional term derived from the PSPG method.    

\subsection{Newton's Method}\label{femnewton}
It is obvious both (\ref{eq:stabweakfluid}) 
and (\ref{eq:stabweakelas}) are nonlinear systems of stabilized 
finite element equations. To solve these nonlinear 
equations, we employ Newton's method, that requires Jacobian 
computations of the nonlinear equations at the current state variables. 
The equations are linearized using standard concept of 
directional derivatives in continuum mechanics; 
see, e.g., \cite{JB08:00, GAH00:00, RW97:00}. 
Since both (\ref{eq:stabweakfluid}) and (\ref{eq:stabweakelas}) 
have similar forms, let $w_k$ denote the fluid velocity or 
structure displacement, 
and $p_k$ be the fluid or structure pressure at 
$k$th iteration. For simplicity of notations, we 
use unified notations $R^1((w_k, p_k), v)$ and $R^2((w_k, p_k), q)$ to 
represent the residuals at the $k$th iteration ($k\ge 0$) in 
the fluid sub-problem (\ref{eq:stabweakfluid}) or the 
residuals in the structure sub-problem (\ref{eq:stabweakelas}). 
At each Newton iteration, we need to find corrections $\delta_w$ 
and $\delta_p$ of the following linearized equations 
as the mixed formulation below:
\begin{subequations}\label{eq:linsys}
\begin{eqnarray}
  \frac{\partial R^1((w_k, p_k), v)}{\partial w}\cdot\delta_w+
  \frac{\partial R^1((w_k, p_k), v)}{\partial p}\cdot\delta_p&\nonumber\\&=
  -R^1((w_k, p_k), v),\\
  \frac{\partial R^2((w_k, p_k), q)}{\partial w}\cdot\delta_w+
  \frac{\partial R^2((w_k, p_k), q)}{\partial p}\cdot\delta_p&\nonumber\\&=
  -R^2((w_k, p_k), q),
  \end{eqnarray}
\end{subequations}
for all $v\in V_{f,h}^0/V_{s,h}$ and $q \in Q_{f,h}/Q_{s,h}$. 
The solution is then updated by $w_{k+1}=w_k+\delta_w$ 
and $p_{k+1}=p_k+\delta_p$. 

To express (\ref{eq:linsys}) as a system of 
linear algebraic equations, we introduce 
two finite elements basis $\{\phi^i\}$ and $\{\varphi^i\}$, $i=1,...,m$, 
for the fluid velocity/structure displacement, and fluid/structure 
pressure, respectively. We choose the 
standard hat function on tetrahedral elements 
as the basis function for both the fluid velocity/structure displacement 
and the fluid/structure pressure. 
The finite element functions $w_{k}$ and $p_{k}$ at Newton iteration $k$ 
can be represented as linear combinations of the basis 
\begin{equation}\label{eq:basis0}
  w_k=\displaystyle\sum_{i=1}^{3m}w_k^i\phi^i \textup{ \quad and\quad  } 
  p_k=\displaystyle\sum_{i=1}^{m}p_k^i\varphi^i,
\end{equation}
where $w_k^i$ and $p_k^i$ are the degrees of freedom, and $m$ is the 
number of vertices of the mesh. In an analogous manner, the corrections 
$\delta_w$ and $\delta_p$ at Newton iteration $k$ are then 
expressed as linear combinations of the basis
\begin{equation}\label{eq:basis1}
  \delta_w=\displaystyle\sum_{i=1}^{3m}\Delta w_k^i\phi^i \textup{ \quad and\quad  } 
  \delta_p=\displaystyle\sum_{i=1}^{m}\Delta p_k^i\varphi^i,
\end{equation}
where $\Delta w_k=(\Delta w_k^1,..., \Delta w_k^{3m})^T$ 
and $\Delta p_k=(\Delta p_k^1,..., \Delta p_k^{m})^T$ are 
the coefficient vectors of Newton increments. 
The global linear algebraic 
system at $k$th iteration is then written in the following block structure:
\begin{equation}\label{eq:globsys}
  K_k\left[
    \begin{array}{c}
      \Delta w_k\\
      \Delta p_k
    \end{array}
  \right]=
  \left[
  \begin{array}{cc}
    A&B_1^T\\
    B_2 &-C
  \end{array}
  \right]
  \left[
    \begin{array}{c}
      \Delta w_k\\
      \Delta p_k
    \end{array}
  \right]
  =
  \left[
    \begin{array}{c}
      r_{1,k}\\
      r_{2,k}
    \end{array}
  \right],
\end{equation}
where $K_k:=K(w_k, p_k)$ is the computed Jacobian matrix containing 
four block matrices $A$, $B_1$, $B_2$ and $C$. For the 
fluid sub-problem, we have $A\neq A^T$, $B_1\neq B_2$, and 
$C=C^T$. For the structure sub-problem, we have $A=A^T$, 
$B_1\neq B_2$, $C=C^T$. In any case, we have to deal 
with an unsymmetric algebraic equation with $K_k\neq K_k^T$. 

Then, Newton's method consists of the following steps: given 
an initial guess $w_0$ and $p_0$, for $k\ge 0$, 
\begin{itemize}
\item[1.] assemble and solve the linearized system (\ref{eq:globsys}) up to 
  a relative residual error reduction factor $\varepsilon_1$, 
\item[2.] update the solutions by $\tilde{w}_{k+1}=w_k+\delta_w$ 
  and $\tilde{p}_{k+1}=p_k+\delta_p$, and go to step 1. 
  until 
  \begin{equation}\label{newstop}
    \|(\tilde{w}_{k+1}-w_k, \tilde{p}_{k+1}-p_k)\|_{L_2}
    \leq\varepsilon_2\|(\tilde{w}_1-w_0, \tilde{p}_1-p_0)\|_{L_2}
  \end{equation}
  is fulfilled. 
\end{itemize}
The Newton iteration is always terminated if the 
stopping criterion (\ref{newstop}) is fulfilled with $\varepsilon_2=10^{-8}$. 
In the first part of our numerical studies presented in Section \ref{sec:ne}, 
we choose $\varepsilon_1:=10^{-8}$ as stopping criterion for 
the solvers of the linearized problems, whereas in the second part, i.e. 
in Subsection \ref{sec:adpt_err}, 
we use an adaptive error control on inner and outer 
iterations in order to improve the overall efficiency. 

More precisely, when employing 
the AMG method as the inner linear solver, with the matrix-vector notations, 
the above nonlinear stopping criterion (\ref{newstop}) can be 
reformulated as follows. 

Let $v_{k}=\left[w_k^T, p_k^T\right]^T$ 
denote the solution at the $k$th Newton iteration for the 
nonlinear system arising from the finite element 
discretization of (\ref{eq:stabweakfluid}) or (\ref{eq:stabweakelas}):
\begin{equation}
  F(v)=0,
\end{equation}
where $v=[w^T, p^T]^T$ denotes the exact solution we are aiming to find. Then 
under appropriate assumption, see, e.g., \cite{PD05:00}, 
we have the following quadratic convergence 
rate of the Newton iteration
\begin{equation}
  \|v-v_{k+1}\|\leq C\|v-v_k\|^2,
\end{equation}
where the norm $\|\cdot\|$ has to be chosen accordingly, and the 
constant $C>0$ may depend on the mesh size on the discrete level and the 
choice of norms. The correction step of the 
Newton iteration is given by 
\begin{equation}
  v_{k+1} = v_k + \Delta v_k
\end{equation}
with the exact correction 
\begin{equation}
  \Delta v_k:=(F^{'}(v_k))^{-1}(-F(v_k))=K_k^{-1}r_k.
\end{equation}

Let $\mathcal{M}_k$ denote the AMG iteration 
matrix at the $k$th Newton iteration. By applying AMG iterations to 
the linearized system (\ref{eq:globsys}) with an initial guess $0$, 
we then have the following relation
\begin{equation}\label{discorr}
  \Delta\tilde{v}_k = {\mathcal C}_k^{-1}r_k
\end{equation}
with the approximated solution 
$\Delta\tilde{v}_k:=\left[\Delta \tilde{w}_k^T, 
  \Delta \tilde{p}_k^T\right]^T$. The iteration matrix ${\mathcal C}_k^{-1}$ 
of the AMG method is given by  
\begin{equation}
  {\mathcal C}_k^{-1} = (I-{\mathcal M}_k^{l(k)})K_k^{-1}
\end{equation} 
with the function $l(k)$ indicating the number of AMG iterations 
for solving the linearized system. 
Since we have only approximate solutions $\tilde{v}^{k+1}$
by applying AMG iterations, we obtain 
\begin{equation}
  \tilde{v}_{k+1}=v_k + \Delta\tilde{v}_k.
\end{equation}

Now we want to estimate the difference $v-\tilde{v}_{k+1}$ 
in the proper norm in order to control the relation 
between the inner and outer accuracy. This can be done 
in the following way
\begin{equation}\label{proof1}
  \begin{aligned}
    \|v-\tilde{v}_{k+1}\|&=\|v-v_{k+1}+v_{k+1}-\tilde{v}_{k+1}\|\\
    &\leq\|v-v_{k+1}\|+\|v_{k+1}-\tilde{v}_{k+1}\|\\
    &\leq C\|v-v_{k}\|^2+\|\Delta v_k-\Delta \tilde{v}_k\|\\
    &= C\|v-v_{k}\|^2+\|K_k^{-1}r_k-{\mathcal C}_k^{-1}r_k\|\\
    &= C\|v-v_{k}\|^2+\|{\mathcal M}_k^{l(k)}K_k^{-1}r_k\|\\
    &\leq C\|v-v_{k}\|^2+\|{\mathcal M}_k\|^{l(k)}\|K_k^{-1}r_k\|
  \end{aligned}
\end{equation}
Now using the relation that $K_k^{-1}r_k=\Delta v_k=v_{k+1}-v_k$ and 
the quadratic convergence of the Newton iteration, we have 
\begin{equation}\label{proof2}
  \begin{aligned}
  \|K_k^{-1}r_k\|=\|v_{k+1}-v_k\|&=\|v-v_{k+1}\|+ \|v-v_k\|\\
  &\leq C\|v-v_k\|^2+ \|v-v_k\|.
  \end{aligned}
\end{equation}
Combing (\ref{proof1}) and (\ref{proof2}), we arrive at 
\begin{equation}
  \begin{aligned}
    \|v-\tilde{v}_{k+1}\|&\leq C\|v-v_{k}\|^2+\|{\mathcal M}_k\|^{l(k)}
    \left(C\|v-v_k\|^2+ \|v-v_k\|\right)\\
    &\leq 2C\|v-v_k\|^2 + \|{\mathcal M}_k\|^{l(k)}\|v-v_k\|.
  \end{aligned}
\end{equation}
Now it is easy to see, that in order to avoid the deterioration of the 
quadratic convergence rate of the Newton iteration, the 
following relation has to be satisfied:
\begin{equation}
  \|{\mathcal M}_k\|^{l(k)}\leq c\|v-v_k\|
\end{equation}
with a constant $c$. 

However, the quantities in the above estimates are not 
computable since we have no exact solution $v$ available. Thus, 
in practice, for the outer iteration, 
we use the $L_2$ norm of the increment in (\ref{newstop}) that can 
be reformulated as 
\begin{equation}
  \begin{aligned}
    \|(\tilde{w}_{k+1}-w_k, \tilde{p}_{k+1}-p_k)\|_{L_2}
    &=\left\| \Delta \tilde{v}_k\right\|_M
    =\left\|\mathcal{C}_k^{-1}r_k\right\|_M\\
    &=\left(M\mathcal{C}_k^{-1}r_k, 
    {\mathcal C}_k^{-1}r_k\right)_{l_2}^{1/2}\\
    &= \left({\mathcal C}_k^{-T}M\mathcal{C}_k^{-1}r_k, 
    r_k\right)_{l_2}^{1/2}\\
    &=\left\|r_k\right\|_{\mathcal C^{-T} M\mathcal{C}_k^{-1}}\\
    &\leq\varepsilon_2\left\|r_0\right\|_{\mathcal C_0^{-T} M\mathcal{C}_0^{-1}},
  \end{aligned}
\end{equation}
where $M=\text{diag}[M_1, M_2]$ 
denotes the mass matrix with the velocity/displacement 
mass matrix $M_1$ and pressure mass matrix $M_2$ on the diagonal. 
From this point of view, (\ref{newstop}) is nothing but 
the relative residual error 
reduction in the ${\mathcal C}_k^{-T}M{\mathcal C}_k^{-1}$-energy norm 
with the symmetric and positive definite matrix 
${\mathcal C}_k^{-T}M{\mathcal C}_k^{-1}$. 

For the inner iteration, we use the relative error reduction 
of the AMG iterations. On the one hand, we have
\begin{equation}
  \begin{aligned}
    \|\Delta v_k-\Delta\tilde{v}_k\|_{K_k^TMK_k}&=\|K_k^{-1}r_k-{\mathcal C}_k^{-1}r_k\|_{K_k^TMK_k}\\
    &=\|K_k^{-1}r_k-(I-{\mathcal M}_k^{l(k)})K_k^{-1}r_k\|_{K_k^TMK_k}\\
    &=\|{\mathcal M}_k^{l(k)}K_k^{-1}r_k\|_{K_k^TMK_k}\\
    &\leq\|{\mathcal M}_k\|_{K_k^TMK_k}^{l(k)}\|K_k^{-1}r_k\|_{K_k^TMK_k}\\
    &=\|{\mathcal M}_k\|_{K_k^TMK_k}^{l(k)}\left(MK_KK_k^{-1}r_k, K_kK_k^{-1}r_k\right)^{1/2}\\
    &=\|{\mathcal M}_k\|_{K_k^TMK_k}^{l(k)}\|r_k\|_M.
  \end{aligned}
\end{equation} 
On the other hand, we have 
\begin{equation}
  \begin{aligned}
    \|\Delta v_k-\Delta\tilde{v}_k\|_{K_k^TMK_k}
    &=\left(K_k^TMK_k(\Delta v_k-\Delta\tilde{v}_k), (\Delta v_k-\Delta\tilde{v}_k)\right)^{1/2}\\
    &=\left(MK_k(\Delta v_k-\Delta\tilde{v}_k), K_k(\Delta v_k-\Delta\tilde{v}_k)\right)^{1/2}\\
    &=(Mr_{AMG}, r_{AMG})^{1/2}\\
    &=\|r_{AMG}\|_M.
  \end{aligned}
\end{equation}
Thus, for the relative residual norm of the linearized system, we have 
the following estimate
\begin{equation}
  \|r_{AMG}\|_M\leq \|{\mathcal M}_k\|_{K_k^TMK_k}^{l(k)}\|r_k\|_M,
\end{equation}
where $r_{AMG}:=r_k-K_k\Delta\tilde{v}_k$ denotes 
the residual of the linearized system in the AMG iterations. 
Since on the discrete level, all the norms are equivalent (up to some 
scaling), we have chosen the error reduction factor of the inner 
iteration such that 
$$\|{\mathcal M}_k^{l(k)}\|_{K_k^TMK_k}\leq\|{\mathcal M}_k\|_{K_k^TMK_k}^{l(k)}
\leq\varepsilon_{2,k-1}$$ is fulfilled, 
where $\varepsilon_{2,k-1}$ indicates the 
error reduction factor at the previous $(k-1)$th Newton iteration. 
The relative residual error reduction of 
the inner iterations will then have the following form
\begin{equation}
  \|r_{AMG}\|_M\leq \varepsilon_{2,k-1}^2=:\varepsilon_{1,k}.
\end{equation}
For more details see \cite{UL13:00,UJ:89,UJ:91,UH:02}. 

In our numerical examples in Section \ref{sec:ne}, we have always 
$\|(w_1-w_0, p_1-p_0)\|_{L_2}\approx 100$. Thus we simply require to 
stop the iteration if 
$\|(w_{k+1}-w_k, p_{k+1}-p_k)\|_{L_2}\leq\varepsilon_2=100\cdot\tilde{\varepsilon}_2$.

\section{Solution methods}\label{sec:sm}
\subsection{A Dirichlet-Neumann FSI iteration}
To solve the nonlinear FSI coupling in this work, 
we use a classical Dirichlet-Neumann (DN) FSI iteration 
(see, e.g., \cite{UK08:00}), 
that is very easy and efficient to implement. Each DN FSI iteration 
consists of the following three steps:
\begin{itemize}
\item[1.] solve the harmonic extension problem with prescribed 
  fluid domain displacement on $\Gamma^0$ from the solution 
  of the structure sub-problem, compute the current 
  fluid domain $\Omega_f^{t^n}$ and 
  the fluid domain velocity $w^{t^n}$, 
\item[2.] with the prescribed fluid velocity condition on $\Gamma^{t^n}$ 
  and the computed fluid sub-domain velocity $w^{t^n}$, 
  and in the computed current fluid domain $\Omega_f^{t^n}$, 
  solve the nonlinear fluid sub-problem (\ref{eq:stabweakfluid}) 
  using Newton's method, 
\item[3.] compute the coupling force on $\Gamma^0$ 
  from the solutions of the 
  fluid sub-problem, solve the structure sub-problem (\ref{eq:stabweakelas}) 
  using Newton's method with the computed force. 
\end{itemize}
Repeat the above three steps until the stopping criterion 
$|r_{\Gamma,k}^n|/\sqrt{n}<\varepsilon$ is reached; see \cite{UK08:00}. 
Here we set $\varepsilon=10^{-8}$ in the numerical experiments, 
$n$ is the number of vertices on the interface $\Gamma^0$ and 
$r_{\Gamma,k}^n$ the interface residual introduced by 
$r_{\Gamma,k}^n = \tilde{d}_{\Gamma, k}^n-d_{\Gamma, k-1}^n$, where 
$d_{\Gamma, k-1}$ denotes the previously 
computed interface displacement on $\Gamma^0$, and $\tilde{d}_{\Gamma, k}^n$ 
the currently computed structure displacement restricted to $\Gamma^0$. 
\subsection{Aitken relaxation}
As is well known, the above simple DN FSI iteration may lead to very 
slow convergence due to the added-mass effect 
(see, e.g., \cite{PC05:00}).  Using the 
Aitken relaxation we are able to 
guarantee and accelerate the convergence of the DN iteration. The 
relaxation step after each DN cycle reads: 
\begin{equation}\label{eq:acceit}
  d_{\Gamma, k}^n = \omega_{k-1}\tilde{d}_{\Gamma, k}^n+(1-\omega_{k-1})d_{\Gamma, k-1}^n,
\end{equation}
where the acceleration parameter $\omega_{k}$ is defined 
by a recursion (see \cite{UK08:00})
\begin{equation}\label{eq:dnomega}
  \omega_k = -\omega_{k-1}\frac{(r^n_{\Gamma, k})^T(r^n_{\Gamma, k+1}-r^n_{\Gamma, k})}
  {|r^n_{\Gamma, k+1}-r^n_{\Gamma, k}|^2}.
\end{equation}
The sequence $d_{\Gamma, k}^n$, $k=1, 2, ...$ converges to a point $d_{\Gamma}^n$ 
that is a fixed-point of the following FSI interface equation:
\begin{equation}\label{eq:fixedp}
  d_{\Gamma}^n=\mathcal{S}^{-1}(\mathcal{F}(d_{\Gamma}^n)), 
\end{equation}
where $\mathcal{F}$ denotes the fluid 
Dirichlet-to-Neumann mapping from the fluid velocity 
to the fluid force on $\Gamma^0$ by 
solving (\ref{eq:stabweakfluid}), and $\mathcal{S}^{-1}$ the 
Neumann-to-Dirichlet mapping from the structure forces to 
the structure displacement on $\Gamma^0$ by solving 
(\ref{eq:stabweakelas}). Using the above mapping notations, 
the fixed-point FSI iteration with relaxation is given by
\begin{equation}\label{eq:fixepr}
  d_{\Gamma, k}^n=\omega_{k-1}\mathcal{S}^{-1}(\mathcal{F}(d_{\Gamma, k-1}^n))
  +(1-\omega_{k-1})d_{\Gamma, k-1}^n.
\end{equation} 
For more details, we refer to \cite{UK08:00}.
\subsection{A short remark on other FSI iterations}
In case of a semi-implicit coupling, where we use a first order 
extrapolation to approximate the current fluid domain and semi-implicit 
treatment for the fluid convection term, at each time step 
the Robin-Neumann (RN) preconditioned 
Krylov subspace method has been used in \cite{SB08:00, SB09:00, YH11:00}. 
This method can be viewed as a generalization of the DN FSI iteration 
(in the semi-implicit coupling setting), that contains 
a more sophisticated preconditioner with weighted 
contributions from both fluid and structure sides. Proper 
choices of involved weighting parameters, combined with 
Krylov subspace acceleration, lead to very fast convergence 
of the FSI iterations. In a fully implicit scheme, we may run a two-level 
approach, that includes nested iterations. 
The outer iteration handles the geometrical and fluid 
convective nonlinearity using first order extrapolation, 
and the inner iteration runs the 
RN preconditioned Kyrlov subspace method. In this work, we stick to the 
above simple and efficient DN FSI iterations with Aitken relaxation. 
\subsection{Preconditioned Kyrlov subspace methods}
The remaining cost in the DN FSI iteration is to solve 
the fluid and structure sub-problems. At each DN iteration, 
we need to 
solve the highly nonlinear fluid and structure 
sub-problems using Newton's method. 
Each Newton iteration requires an efficient solver 
for the linearized system 
(\ref{eq:globsys}). From now on, for simplicity, 
instead of $K_k$ we use the notation $K$ to indicate the 
linearized Jacobian. We first consider some Krylov subspace methods with 
efficient preconditioners for the linear system of algebraic equations, 
that we use in this work:
\begin{equation}\label{eq:prelinsys}
   P^{-1}K\left[
    \begin{array}{c}
      \Delta w_k\\
      \Delta p_k
    \end{array}
  \right]=
  P^{-1}\left[
    \begin{array}{c}
      r_{1,k}\\
      r_{2,k}
    \end{array}
  \right],
\end{equation} 
where $P$ is a properly chosen preconditioner. 

An efficient preconditioner for the linear system (\ref{eq:globsys}) is 
developed according to the following $LU$ factorization 
of the matrices in the linear 
system of the form (\ref{eq:globsys}) (see, e.g., \cite{PSV08:00, YS03:00}):
\begin{equation}\label{eq:lufacex} 
  \left[\begin{array}{cc}
      I& 0\\
      B_2A^{-1}& I
    \end{array}
    \right]
  \left[\begin{array}{cc}
      A& 0\\
      0&-S
    \end{array}
    \right]
  \left[\begin{array}{cc}
      I& A^{-1}B^T_1\\
      0 & I
    \end{array}
    \right]=:LDU,
\end{equation}
where $S=B_2A^{-1}B_1^T+C$ is the (negative) Schur complement. The 
inverses of $A$ and $S$ are very expensive in general. Thus we construct 
efficient preconditioner based on the following modified $LU$ factorization. 
\begin{equation}\label{eq:lufac} 
  \left[\begin{array}{cc}
      I& 0\\
      B_2\tilde{A}^{-1}& I
    \end{array}
    \right]
  \left[\begin{array}{cc}
      \tilde{A}& 0\\
      0&-\tilde{S}
    \end{array}
    \right]
  \left[\begin{array}{cc}
      I& \tilde{A}^{-1}B_1^T\\
      0 & I
    \end{array}
    \right]=:\tilde{L}\tilde{D}\tilde{U},
\end{equation}
where $\tilde{A}^{-1}$ and 
$\tilde{S}^{-1}$ are proper approximations for $A^{-1}$ and ${S}^{-1}$, 
that are easy to invert.
\subsubsection{A preconditioner for the fluid sub-problem}
To solve the fluid sub-problem, we use a Generalized 
Conjugate Residual (GCR) method (see, e.g., \cite{YS03:00}) 
with a preconditioner $P_R=\tilde{D}\tilde{U}$. The inverse of the preconditioner 
is given by 
\begin{equation}\label{eq:rpre}
  P_R^{-1}=
  \left[\begin{array}{cc}
      \tilde{A}^{-1}& \tilde{A}^{-1}B_1^T\tilde{S}^{-1}\\
      0& -\tilde{S}^{-1}
    \end{array}
  \right],
\end{equation}
where $\tilde{A}^{-1}$ is realized by applying AMG cycles 
(see, e.g., \cite{JW87:00, FK98:00}), and 
$\tilde{S}^{-1}$ is computed by the operator splitting technique 
(see, e.g., \cite{ST99:00, AW08:00}). 
The method starts with the splitting of the velocity matrix $A$
\begin{equation}\label{eq:split}
  A=\frac{\rho_f}{\Delta t}M+\mu D + \rho_f C,
\end{equation}
where $M$ stands for velocity mass matrix, $D$ the diffusion matrix, $C$ the 
convection matrix. Using the properties of corresponding operators, 
the Schur complement is then approximated by the 
following addition 
\begin{equation}\label{eq:rpreschur}
  \begin{aligned}
  (B_2A^{-1}B_1^T)^{-1}\approx&\frac{\rho_f}{\Delta t}(B_2M^{-1}B_1^T)^{-1}
  +\mu(B_2D^{-1}B_1^T)^{-1}
  +\rho_f(B_2C^{-1}B_1^T)^{-1}\\
  \approx&\frac{\rho_f}{\Delta t}D_p^{-1}+\mu M_p^{-1}
  +\rho_fM_p^{-1}C_pD_p^{-1}=:\tilde{S}^{-1},
\end{aligned}
\end{equation} 
where $D_p$ and $C_p$ are the stiffness matrices associated with the 
finite element discretization of the Laplacian operator and the scalar 
convection operator in the pressure space, respectively, 
$M_P$ the mass matrix in the pressure space. In actual calculations, 
$D_p$ is inverted by AMG method, and $M_p$ 
is replaced by $\textup{diag}(M_p)$, i.e., the diagonal of $M_p$. 
These operations are relatively cheap to realize. In \cite{EM10:00}, 
this method was also used to solve the linearized NS 
equations using the $P_2-P_1$ mixed velocity and pressure formulation. 

It is easy to see, that the operation of $P_R^{-1}$ applied to a vector 
consists of 
two steps: solve the pressure Schur complement 
equation with $\tilde{S}$, and then solve the velocity 
equation with $\tilde{A}$. 
\subsubsection{A preconditioner for the structure sub-problem}
To solve the structure sub-problem, we use the 
Biconjugate Gradient Stabilized (BiCGStab) method (see \cite{HV92:00}) 
with a preconditioner $P_L=\tilde{L}\tilde{D}$. 
The inverse of the preconditioner 
is given by 
\begin{equation}\label{eq:lpre}
  P_L^{-1}=
  \left[\begin{array}{cc}
      \tilde{A}^{-1}& 0\\
      \tilde{S}^{-1}B_2\tilde{A}^{-1}& -\tilde{S}^{-1}
    \end{array}
  \right],
\end{equation}
where $\tilde{A}^{-1}$ is performed by applying AMG cycles 
to $A$, and the Schur complement is approximated by 
\begin{equation}\label{eq:lpreschur}
  S\approx\tilde{S}=(\frac{1}{\theta}+\frac{1}{\kappa})\textup{diag}(M_p),
\end{equation}
where $M_p$ is the mass matrix in the structure pressure space. 
This preconditioner has demonstrated the robustness with 
respect to near-incompressibility. 
In our numerical experiments for both hyperslastic models, 
we use $\theta=O(c_{10})$. A good choice of $\theta$ is 
usually adjusted from the numerical tests on a coarse mesh. 
Once chosen, the value is held fixed for finer meshes. 
See \cite{ME10:00} for the application of this method 
using $P_2-P_1$ mixed displacement 
and pressure formulation in 
the modified Mooney-Rivlin hyperelastic model.  

It is obvious to see, that the 
operation of $P_L^{-1}$ applied to a vector 
consists of two steps: solve the velocity 
equation with $\tilde{A}$ 
and then solve the pressure equation 
with the approximated Schur complement 
$\tilde{S}$. 
\subsection{AMG methods for the mixed problems}
To improve performance of the linear iterative solvers, a class of 
(nearly) optimal AMG methods for the $P_1-P_1$ 
mixed formulation are considered. 
In \cite{WM04:00, WM06:00, Yang11:00, YH11:00}, the robust coarsening 
strategy is discussed for 
stabilized $P_1-P_1$ discretization of the 
\textit{Oseen} equations and nearly incompressible linear elastic equations, 
that guarantees the discrete $\inf-\sup$ conditions on all coarse levels in 
an algebraic manner. 
In this work, we will focus on proper choices of smoothers in the AMG 
methods for both fluid and structure sub-problems. As observed 
there exists no \textit{black-box} solution, that can be applied 
to both sub-problems. In the AMG method, we need 
both pre- and post-smoothing steps, that are applied to all levels. For 
simplicity of presentation, we only discuss a smoothing step on one particular 
level. The smoothing step on a coarse level is 
a straightforward application to the corresponding linear system 
obtained by the stabilized Galerkin projection.
\subsubsection{The {\it Braess-Sarazin} Smoother}
The {\it Braess-Sarazin} consists of three steps of 
the {\it{inexact symmetric Uzawa algorithm}}:
\begin{subequations}\label{eq:rismsteps}
  \begin{eqnarray}
    \Delta w_k&=&\Delta w_k+\tilde{A}^{-1}(r_{1, k}-A\Delta w_k-B_1^T\Delta p_k),\\
    \Delta p_k&=&\Delta p_k-\underbrace{\tilde{S}^{-1}(r_{2, k}-B_2\Delta w_k+C\Delta p_k)}_{:=r_{p, k}},\\
    \Delta w_k&=&\Delta w_k +\tilde{A}^{-1}B_1^Tr_{p,k},
  \end{eqnarray}
\end{subequations}
which corresponds to a Richardson iteration applied to 
(\ref{eq:globsys}):
\begin{equation}\label{eq:bsglobsys}
  \left[
    \begin{array}{c}
      \Delta w_k\\
      \Delta p_k
    \end{array}
  \right]=
  \left[
    \begin{array}{c}
      \Delta w_k\\
      \Delta p_k
    \end{array}
  \right]
  +P_F^{-1}
  \left(
    \left[
      \begin{array}{c}
        r_{1,k}\\
        r_{2,k}
      \end{array}
    \right]-
    K\left[
      \begin{array}{c}
        \Delta w_k\\
        \Delta p_k
      \end{array}
    \right]  
  \right)
\end{equation}
with the full preconditioner $P_F$ given by
\begin{equation}\label{eq:fpre}
  P_F=
  \left[\begin{array}{cc}
      \tilde{A}& B_1^T\\
      B_2& B_2\tilde{A}^{-1}B_1^T-\tilde{S}
    \end{array}
  \right].
\end{equation}
An essential issue is how to choose preconditioners 
$\tilde{A}$ and $\tilde{S}$ in order to fulfill the 
smoothing property. We use $\tilde{A}=2D$, 
where $D$ denotes the diagonal of $A$ (a relaxed Jacobi iteration). 
The original {\it Braess-Sarazin} 
smoother (see \cite{Braess97:00}) 
needs to solve the pressure correction equation 
exactly with the Schur-complement $\tilde{S}=B_2\tilde{A}^{-1}B_1+C$. A relaxed 
version (see \cite{WZ00:00}) 
only solves this equation approximately, using e.g., an 
inner AMG method with starting value $0$. 
As a price to pay, in each smoothing step of the outer AMG method, 
we first construct 
a Schur complement $\tilde{S}$ on each level using the Galerkin 
projected matrices $B_2$, $\tilde{A}^{-1}$, $B_1$ and $C$. Then 
on each level of the inner AMG method, we 
need to construct a 
Schur complement using standard Galerkin projection applied 
to the previously constructed Schur complement. 
In \cite{WM04:00, WM06:00, Yang11:00, YH11:00}, the {\it Braess-Sarzin} 
smoother demonstrates robustness and efficiency in the AMG 
methods for solving the \textit{Oseen} and linear elastic equations. 
However, as observed in this work, this smoother fails 
in the AMG method for solving the hyperelastic equations.  
\subsubsection{The multiplicative {\it Vanka} smoother}
We finally consider the multiplicative {\it Vanka} smoother. 
In each smoothing step, 
we need to solve a sequency of local problems on properly chosen 
patches $\mathcal{P}_i$, $i=1,..., n$, 
where $n$ is the number of pressure 
degrees of freedom on that level. 
Each patch $\mathcal{P}_i$ contains a pressure degree of freedom, and 
the connected velocity degrees of freedom (given by the connectivity 
of matrix $B_2$). The local problem on $\mathcal{P}_i$ is 
then constructed by a canonical projection of 
the global matrix $K$ on the patch $\mathcal{P}_i$, that 
has a form similar to the global problem (\ref{eq:globsys}): 
\begin{equation}\label{eq:locsys}
  K^i\left[
    \begin{array}{c}
      \Delta w_k^i\\
      \Delta p_k^i
    \end{array}
  \right]=
  \left[
  \begin{array}{cc}
    A^i&(B_1^T)^i\\
    B_2^i &-C^i
  \end{array}
  \right]
  \left[
    \begin{array}{c}
      \Delta w_k^i\\
      \Delta p_k^i
    \end{array}
  \right]
  =
  \left[
    \begin{array}{c}
      r_{1,k}^i\\
      r_{2,k}^i
    \end{array}
  \right],
\end{equation}
where the current local residual $[(r_{1,k}^i)^T ; (r_{2,k}^i)^T]^T$ is 
assembled by incorporating updated solutions from all previously 
treated local problems, which corresponds to a Gauss-Seidel manner. 
Thanks to small size of each local problem, direct solvers are applicable 
efficiently. In order to achieve smoothing property, 
a relaxation is usually necessary 
in the solution updating step:
\begin{equation}\label{eq:vankasmooth}
  \left[
    \begin{array}{c}
      \Delta w_k\\
      \Delta p_k
    \end{array}
    \right]=
  \left[
    \begin{array}{c}
      \Delta w_k\\
      \Delta p_k
    \end{array}
    \right]+
  \omega P_k^i
  \left[
    \begin{array}{c}
      \Delta w_k^i\\
      \Delta p_k^i
    \end{array}
    \right],
\end{equation}   
where $P_k^i$ is the canonical embedding 
from the local patch $\mathcal{P}_i$ to the global, and 
$\omega\in[0.5, 0.9]$ is the relaxation parameter 
used in our numerical experiments. We fix $\omega=0.78$ and $\omega=0.86$ 
for the structure sub-problems modeled by the modified hyperelastic 
Mooney-Rivlin material and the two-layer thick-walled artery, 
respectively. 
Such a {\it Vanka} smoother and its variants have been used in 
the GMG methods for solving both 
the fluid and hyperelastic problems in 2D or 3D, 
see, e.g., \cite{Vanka86:00, ZW03:00, John01:00, John06:00, TS09:00}. 
However, in the AMG method for solving 
the linearized NS equations, the smoothing 
property deteriorates rapidly as observed in our numerical experiments, 
that leads to an invalidity of 
the multigrid convergence; see also the 
failure report of the AMG method 
for the 3D \textit{Oseen} equations in \cite{MW03:00}. For solving the 
linearized hyperelastic equations used in this work, the {\it Vanka} smoother 
works quite well with sufficiently large number of smoothing steps. 

\section{Numerical experiments}\label{sec:ne}
\subsection{Geometrical and material parameters}
For the hyperelastic model of 
two-layer thick-walled artery from a rabbit, 
the geometry and material parameters are chosen 
according to Fig. 14 in \cite{Holzapfel00:00}. The thickness of 
the media and adventitia is set to $0.26$ mm and $0.13$ mm, 
respectively. The angles between the collagen fibers 
and the circumferential direction in the media 
and adventitia are set to $29.0^{\circ}$ and $62.0^{\circ}$, respectively. 
The radius of the artery is set to $1.43$ mm. The length of the 
artery is set to $18$ mm. 

The shear stress-like material parameters $c_{10}$ 
are set to $3$ kPa and $0.3$ kPa, $k_2$ are set to 
$2.3632$ kPa and $0.562$ kPa, the dimensionless parameters 
$k_2$ are set to $0.8393$ and $0.7112$, for the media and adventitia, 
respectively. The bulk modulus is set to $10^5$ kPa for both two layers, 
which corresponds to (nearly) incompressible material. 
The density on the reference configuration is set to $1.2$ mg/mm$^3$. 

For the modified hyperelastic model of Mooney-Rivlin material, 
we adopt the same geometrical configuration as in the artery. We set 
material parameters $c_{10}$ and $c_{01}$ to $3$ kPa and $0.3$ 
kPa, respectively. The bulk modulus is set to $10^5$ kPa. 
The density on the reference configuration is set to $1.2$ mg/mm$^3$.

For the fluid model, we set the density to $1$ mg/mm$^3$, 
the dynamic viscosity to $0.035$ Poise. 

The fluid and structure are 
at the rest in the beginning. The structure 
body is fixed at two ends. 
We give a pressure pulse (Neumann data) for the fluid 
at the innet $\Gamma^t_{in}$ of the fluid domain: $(0, 0, 1.332)$ kPa for 
time $t\leq 1$ ms, and set $(0, 0, 0)$ kPa after $1$ ms in the case of hyperelastic 
Mooney-Rivlin material and  $(0, 0, 1.332)$ kPa for 
time $t\leq 0.125$ ms, and set $(0, 0, 0)$ kPa after $0.125$ ms 
in the case of two-layer thick-walled artery. On the 
outlet $\Gamma^t_{out}$, a doing-nothing condition 
(0 Neumann boundary conditions) is applied. 
Time step size is set $0.125$ ms for both models. 
\subsection{Coarse and fine meshes}
On the fluid coarse mesh, we have $8120$ tetrahedral 
elements, $2259$ vertices and  $9036$ degrees of freedom. 
On the structure coarse mesh, we have $30472$ 
tetrahedral elements, $6524$ vertices and $26096$ degrees of freedom. 
On the fluid fine mesh, we have $64960$ tetrahedral elements, 
$14249$ vertices, and $56996$ degrees of freedom. 
On the structure fine mesh, we have $243776$ 
tetrahedral elements ($131072$ for the media and $112704$ for the 
adventitia), $46356$ vertices, and 
$185424$ degrees of freedom. We use 
the automatic mesh generator Netgen \cite{JS97:00} to generate the FSI meshes, 
that have conforming mesh grids on the FSI interface and resolve different 
layers of hyperelastic models so that different material parameters on 
layers are assigned accordingly. 

\subsection{Nonlinear DN FSI iterations}
\begin{figure}[htbp]
  \centering
  \begin{subfigure}[h]{0.4\textwidth}
    \centering
    \includegraphics[scale = 0.26]{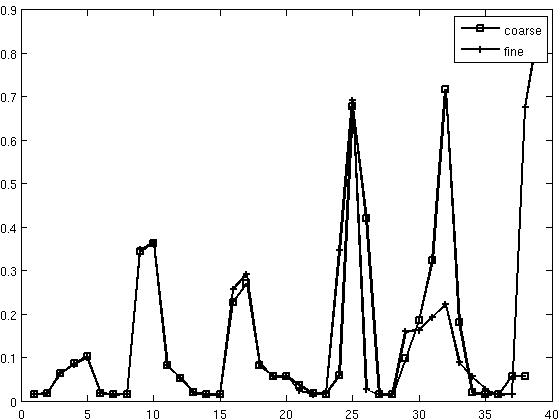}
    \caption{Atiken parameters}\label{mr_par}
  \end{subfigure}
  \hfill
  \begin{subfigure}[h]{0.4\textwidth}
    \centering
    \includegraphics[scale = 0.265]{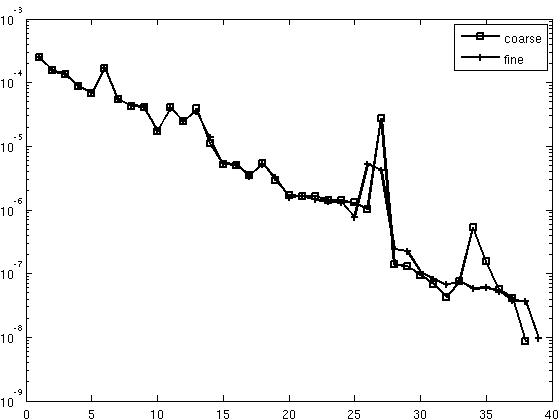}
    \caption{DN convergence history}\label{mr_dnconv}
  \end{subfigure}
  \caption{Nonlinear iterations at the second time step for the Mooney-Rivlin material: Atiken parameters in the nonlinear DN iterations (left) and convergence history of DN FSI iterations (right) on coarse and fine meshes}\label{mr_dn}
\end{figure}

\begin{figure}[htbp]
  \centering{
  \begin{subfigure}[h]{0.4\textwidth}
    \centering
    \includegraphics[scale = 0.26]{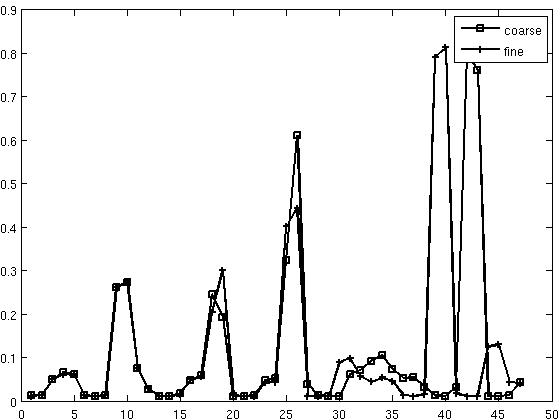}
    \caption{Atiken parameters}\label{ar_par}
  \end{subfigure}
  \hfill
  \begin{subfigure}[h]{0.4\textwidth}
    \centering
    \includegraphics[scale = 0.265]{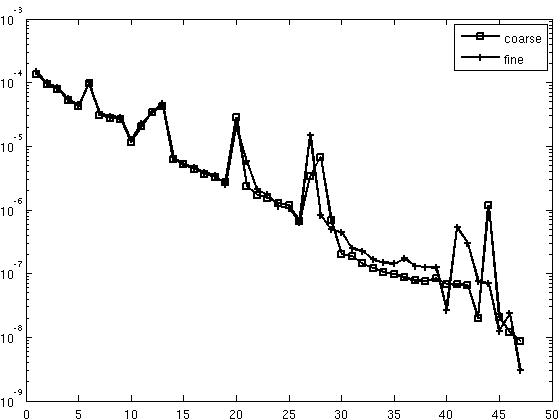}
    \caption{DN convergence history}\label{ar_dnconv}
  \end{subfigure}}
  \caption{Nonlinear iterations at the second time step for the two-layer thick-walled artery: Atiken parameters in the nonlinear DN iterations (left) and convergence history of DN FSI iterations (right) on coarse and fine meshes}\label{ar_dn}
\end{figure}

We first demonstrate the performance of nonlinear DN FSI solvers for 
the Mooney-Rivlin material and two-layer thick-walled 
artery in Fig. \ref{mr_dn} 
and Fig. \ref{ar_dn}, respectively. We depict the value of 
the Atiken relaxation parameter (vertical lines) with respect to 
each DN iteration (horizontal lines) in the left plot, and the 
convergence history (vertical lines) with respect to 
the DN iteration (horizontal lines) in the right plot. Since at 
each time step, the nonlinear DN solver behaves in an analogous manner, we 
only plot the results at the second time step. 

As observed from Fig. \ref{mr_dn} and Fig. \ref{ar_dn}, the DN FSI solver 
for both nonlinear hyperelastic models runs robustly and efficiently. 
The Atiken parameters are updated dynamically according to (\ref{eq:dnomega}) 
at each DN iteration, and are in the range of $(0, 1)$. The iteration 
numbers stay in the same range on both coarse and fine meshes 
for each structure model. 

\begin{figure}[htbp]
  \centering{
    \includegraphics[scale = 0.5]{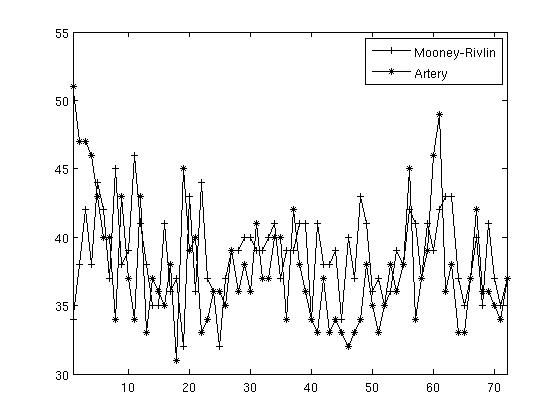}
  }
  \caption{Number of nonlinear DN FSI iterations at different time steps using the model of the modified Mooney-Rivlin material ('-+-') and the model of the two-layer thick-walled artery ('-*-'). The horizontal line indicates the time steps, the vertical line indicates the number of DN iterations.}\label{mrar_dn}
\end{figure}

In Fig. \ref{mrar_dn}, we show the number of nonlinear DN FSI iterations at 
different time steps (up to the $72$th time step), 
using the model of the modified Mooney-Rivlin material 
and the model of the two-layer thick-walled artery. 
As we see from the numerical results, 
the iteration numbers of the DN iteration 
stay in a similar range with different time steps. 

\begin{figure}[htbp]
  \centering
  \begin{subfigure}[h]{0.4\textwidth}
    \centering
    \includegraphics[scale = 0.2]{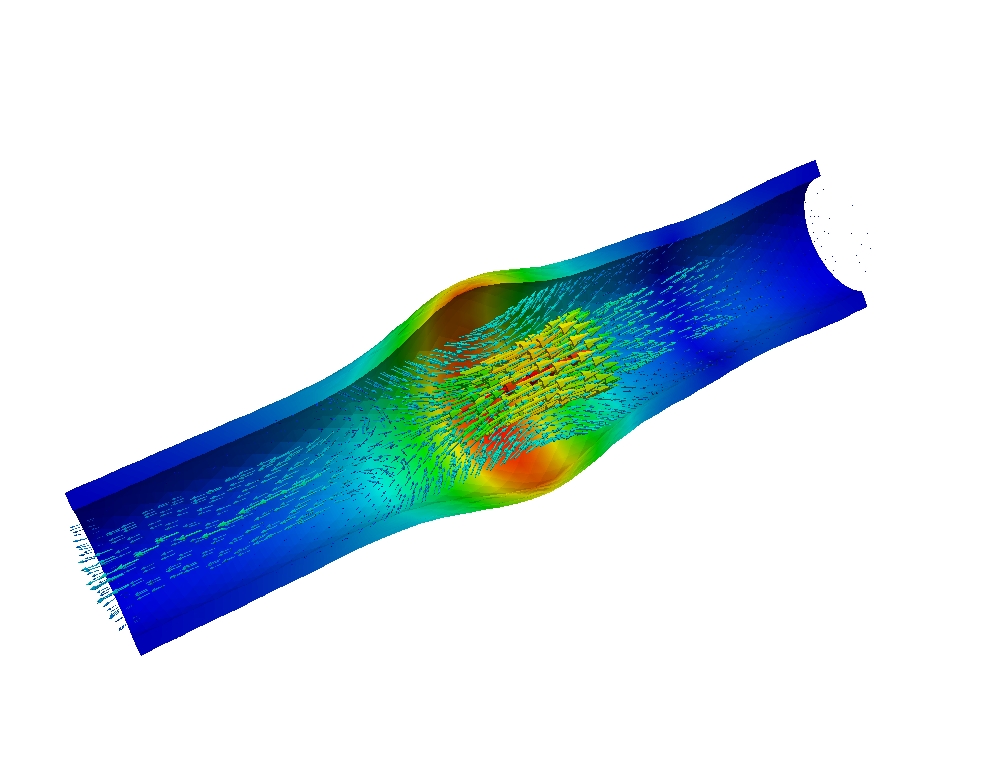}
    \caption{FSI of Mooney-Rivlin}\label{mr_sim}
  \end{subfigure}
  \hfill
  \begin{subfigure}[h]{0.4\textwidth}
    \centering
    \includegraphics[scale = 0.2]{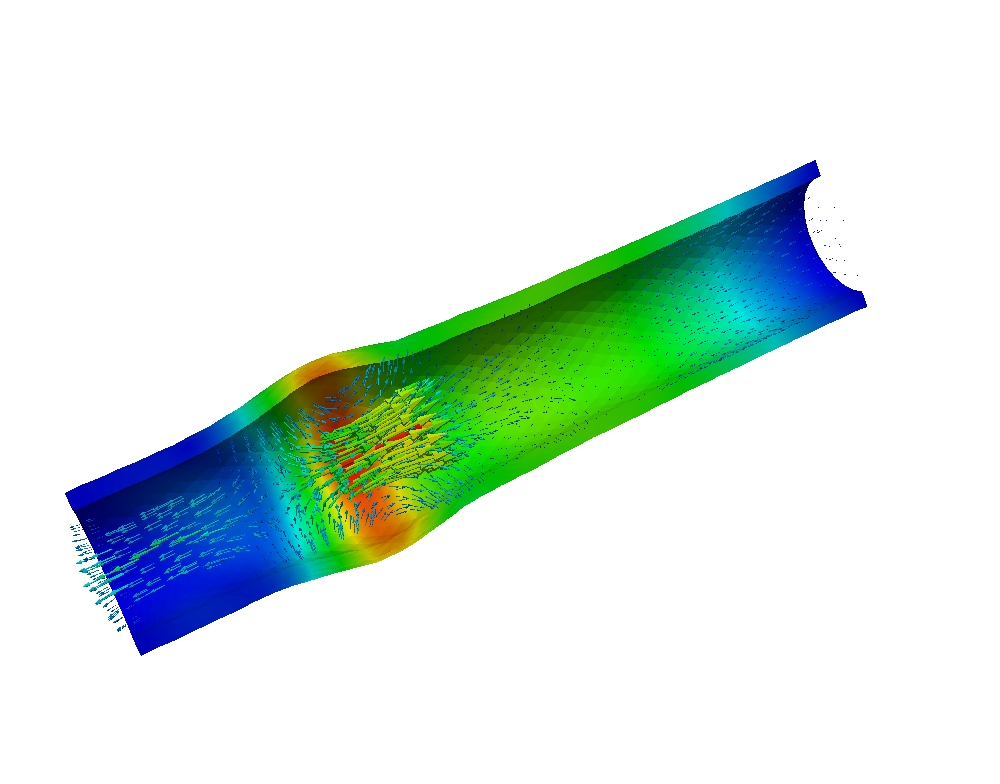}
    \caption{FSI of artery}\label{ar_sim}
  \end{subfigure}
  \caption{FSI simulation using the Mooney-Rivlin material (left) and the two-layer thick-walled artery (right). }\label{fsi_sim}
\end{figure}

For an illustration of the 
FSI simulation results, we visualize the structure deformation, 
the fluid velocity for the Mooney-Rivlin material 
and the two-layer thick-walled artery at time $8$ ms in Fig. \ref{fsi_sim}. 

\subsection{Solution methods for sub-problems}
We first presents numerical results concerning Newton's method combined with 
the Krylov subspace methods: GCR and BiCGStab for the 
fluid and structure sub-problems, respectively. 
After that numerical results of Newton's method combined with 
the AMG methods for both sub-problems are shown. 
For simplicity of presentation, we only present the convergence history 
of Newton's method at the second DN FSI nonlinear iteration at the first 
time step. The iteration numbers of GCR and BiCGStab methods are recorded 
for each Newton iteration. We observe similar behavior of both Newton's method 
and Krylov subspace methods at other time steps. On the left plots 
of the following figures, 
the horizontal lines represent the iteration number of Newton's method, 
the vertical lines 
represent the corresponding errors. On the right plots 
of the following figures, 
the horizontal lines represent the iteration number of Newton's method, 
the vertical lines represent the iteration number of linear solvers needed at 
each Newton iteration. 

\subsubsection{Krylov subspace methods}
\begin{figure}[htbp]
  \centering
  \begin{subfigure}[h]{0.4\textwidth}
    \centering
    \includegraphics[scale = 0.26]{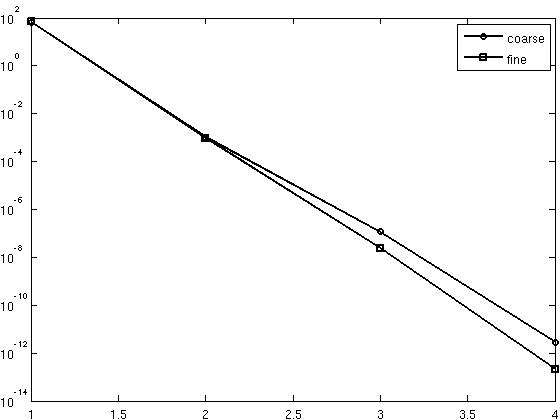}
    \caption{Newton's method}\label{krylovfluid_newton}
  \end{subfigure}
  \hfill
  \begin{subfigure}[h]{0.4\textwidth}
    \centering
    \includegraphics[scale = 0.25]{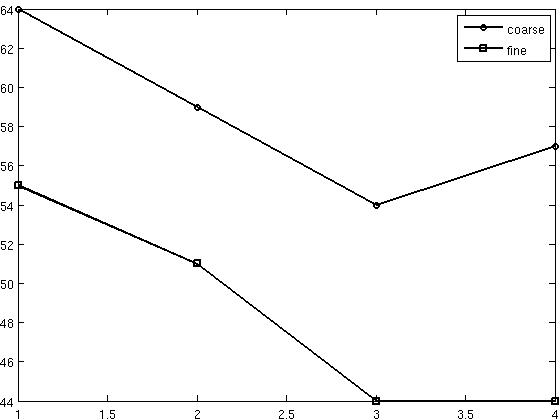}
    \caption{CGR}\label{krylovfluid_gcr}
  \end{subfigure}
  \caption{Convergence history of Newton's method (left) and 
    iteration numbers of GCR solvers (right) for the fluid 
    sub-problem on coarse and fine meshes}\label{krylovfluid}
\end{figure}

The convergence history of Newton's method to solve the NS equations on 
the coarse and fine meshes and the iteration numbers of the GCR method 
to solve the linearized NS equations are recorded in Fig. \ref{krylovfluid}.  

As observed from the numerical results, we obtain 
quadratic convergence rate for Newton's method. The GCR 
method combined with the preconditioner (\ref{eq:rpre}) demonstrates 
the robustness with respect to the mesh refinement, and the material 
parameters. 

\begin{figure}[htbp]
  \centering
  \begin{subfigure}[h]{0.4\textwidth}
    \centering
    \includegraphics[scale = 0.26]{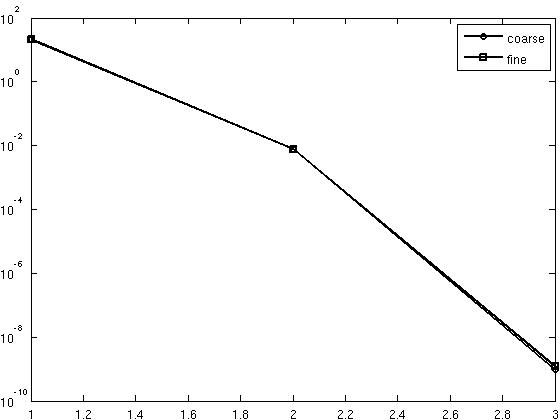}
    \caption{Newton's method}\label{krylovstruc_newton}
  \end{subfigure}
  \hfill
  \begin{subfigure}[h]{0.4\textwidth}
    \centering
    \includegraphics[scale = 0.25]{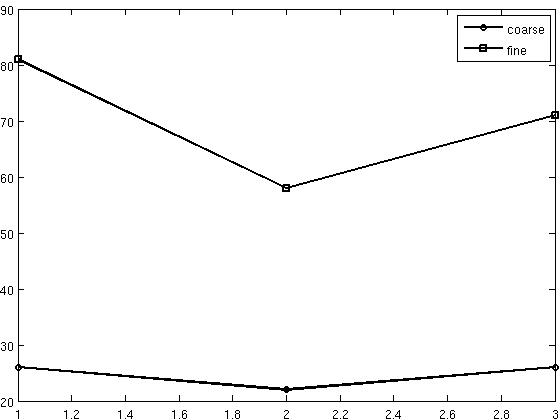}
    \caption{BiCGStab}\label{krylovstruc_bicg}
  \end{subfigure}
  \caption{Convergence history of Newton's method (left) and 
    iteration numbers of  BiCGStab solvers (right) for the structure 
    sub-problem modeled by modified hyperelastic Mooney-Rivlin material 
    on coarse and fine meshes}\label{krylovstruc_mr}
\end{figure}

For the structure sub-problem modeled by the modified hyperelastic 
Mooney-Rivlin material, the convergence history of Newton's method and 
the iteration numbers of the BiCGStab method 
to solve the linearized hyperelastic equations are 
recorded in Fig. \ref{krylovstruc_mr}.  

As observed from the numerical results, we obtain 
quadratic convergence rate for Newton's method. The BiCGStab 
method combined with the preconditioner (\ref{eq:lpre}) shows 
medium number of iteration numbers on both coarse and fine meshes. We 
observe increasing number of iterations with mesh refinement. However, 
the method is robust with respect to the material parameters.  

\begin{figure}[htbp]
  \centering
  \begin{subfigure}[h]{0.4\textwidth}
    \centering
    \includegraphics[scale = 0.26]{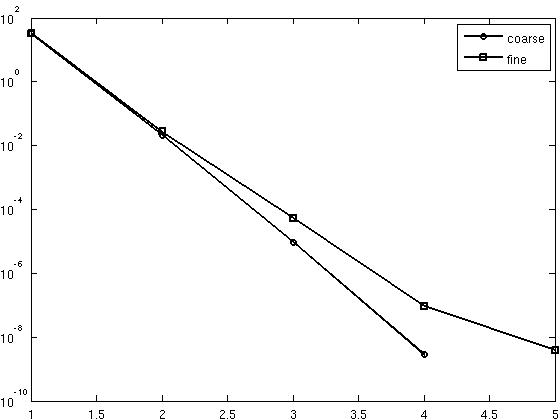}
    \caption{Newton's method}\label{krylovstrucart_newton}
  \end{subfigure}
  \hfill
  \begin{subfigure}[h]{0.4\textwidth}
    \centering
    \includegraphics[scale = 0.25]{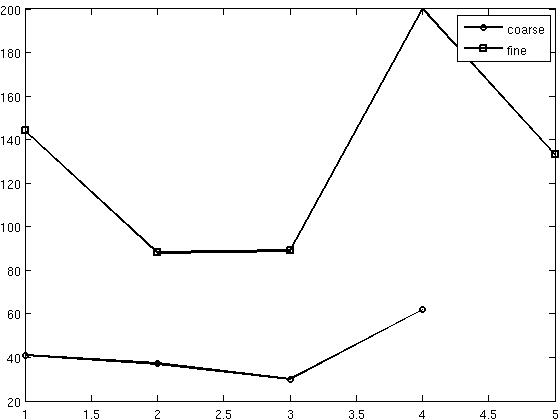}
    \caption{BiCGStab}\label{krylovstrucart_bicg}
  \end{subfigure}
  \caption{Convergence history of Newton's method (left) and 
    iteration numbers of  BiCGStab solvers (right) for the structure 
    sub-problem modeled by the two-layer thick-walled artery 
    on coarse and fine meshes}\label{krylovstruc_artery}
\end{figure}

For the structure sub-problem modeled by the two-layer thick-walled artery, 
the convergence history of Newton's method and 
the iteration numbers of the BiCGStab method 
to solve the linearized hyperelastic equations are 
recorded in Fig. \ref{krylovstruc_artery}.  

As observed from the numerical results, we obtain the quadratic convergence 
of Newton's method on the coarse mesh. On the fine mesh, the quadratic 
convergence degenerates at the fourth step, that is because the BiCGStab method 
at this step 
does not solve the linearized hyperelastic equations accurately. As we see 
from the right plot in Fig. \ref{krylovstruc_artery}, we stop the linear 
solver at $200$ steps if the relative residual error does not reach the factor $10^{-8}$. 
This is cured later on by the more robust AMG method. As in the previous 
hyperelastic model, we observe the medium number of iteration numbers 
on both coarse and fine meshes, and increasing number of iterations with mesh refinement. 

\subsubsection{The AMG methods}
In order to run the AMG methods, we first construct the hierarchy 
of matrices on all levels in pure algebraic way 
(distinct from the GMG methods, where a 
hierarchy of nested meshes are provided). 
We use the coarsening strategy in \cite{FK98:00}. 
On the coarse mesh, we arrive at three and four levels in the 
AMG methods for the fluid and structure sub-problems, respectively. 
For the fluid sub-problem on the coarse mesh, we have $9034$, $1684$ and $252$ degrees of 
freedom on three levels after coarsening. 
For the structure sub-problem, we have $26096$, $4012$, $648$ and $88$ degrees of freedom 
on four levels. 
For the fluid sub-problem on the fine mesh, we have $56996$, 
$9034$, $1684$ and $252$ degrees of freedom on four levels after coarsening. 
For the structure sub-problem, we have $185424$, $26096$, $4012$, 
$648$ and $88$ degrees of freedom on five levels.  
\begin{figure}[htbp]
  \centering
  \begin{subfigure}[h]{0.4\textwidth}
    \centering
    \includegraphics[scale = 0.26]{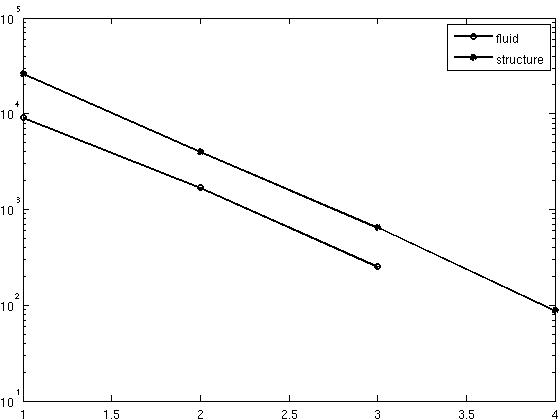}
    \caption{Degrees of freedom}\label{amgcoarsening_coarse}
  \end{subfigure}
  \hfill
  \begin{subfigure}[h]{0.4\textwidth}
    \centering
    \includegraphics[scale = 0.25]{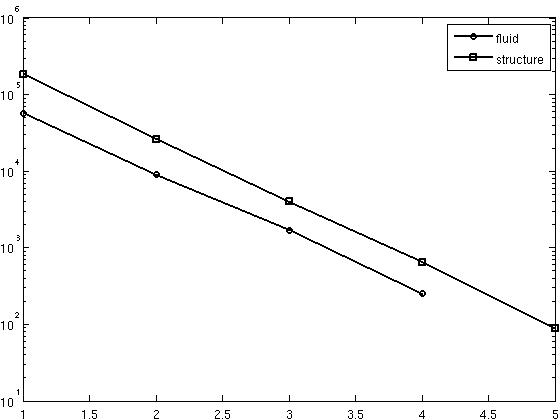}
    \caption{Degrees of freedom}\label{amgcoarsening_fine}
  \end{subfigure}
  \caption{Number of degrees of freedom for the fluid and structure sub-problems 
    on AMG levels using the coarse (left) and fine (right) meshes.}\label{amgcoarsening}
\end{figure}

As we see from Fig. \ref{amgcoarsening}, the number of degrees of freedom 
is reduced by a factor about $8$ from the fine to coarse levels. 

\begin{figure}[htbp]
  \centering
  \begin{subfigure}[h]{0.4\textwidth}
    \centering
    \includegraphics[scale = 0.26]{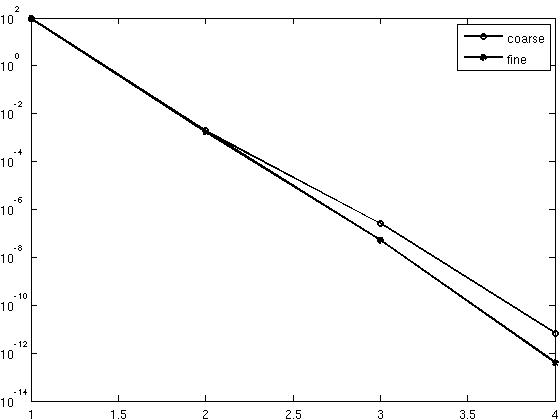}
    \caption{Newton's method}\label{amgfluid_newton}
  \end{subfigure}
  \hfill
  \begin{subfigure}[h]{0.4\textwidth}
    \centering
    \includegraphics[scale = 0.15]{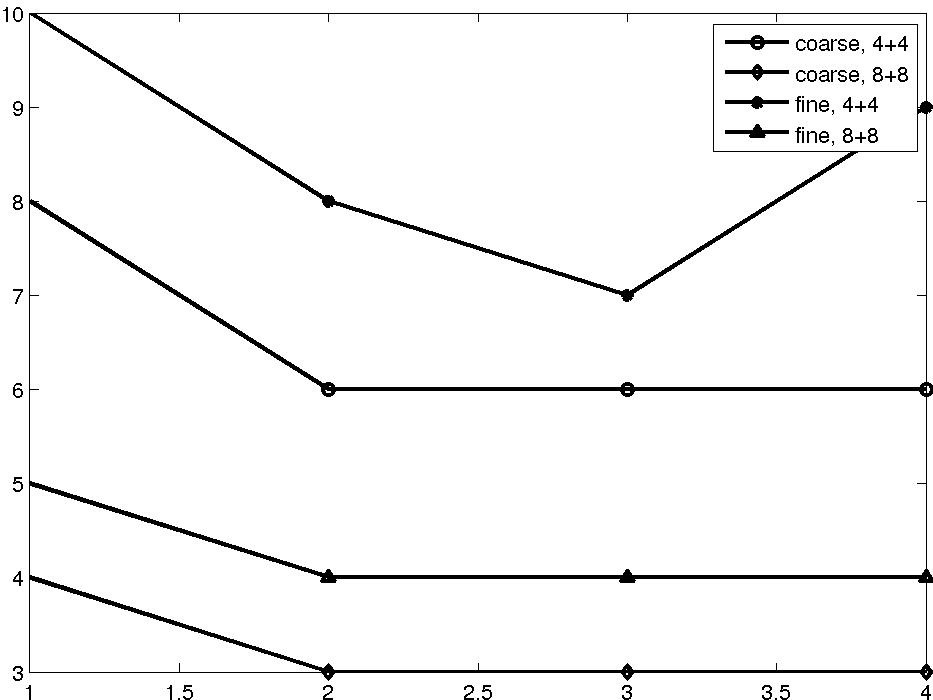}
    \caption{AMG}\label{amgfluidf}
  \end{subfigure}
  \caption{Convergence history of Newton's method (left) and 
    iteration numbers of AMG method (right) for the fluid 
    sub-problem on coarse and fine meshes, using $4$ and $8$ 
    pre- and post-smoothing steps, respectively.}\label{amgfluid}
\end{figure}

For the fluid sub-problem, the convergence of Newton's method and 
the iteration numbers of the AMG method to solve the linearized 
NS equations using different number of {\it Braess-Sarazin} smoothing steps 
are plotted in Fig. \ref{amgfluid}.

To sum from Fig. \ref{amgfluid}, 
a quadratic convergence rate of Newton's method is observed. The number of 
AMG iterations is drastically decreased compared to the GCR method combined 
with the preconditioner. In addition, the iteration numbers are less effected by the 
mesh refinement. By doubling the pre- and post-smoothing steps (from $4$ to $8$), 
we observe the reduced AMG iteration 
numbers by almost a factor $2$. 

\begin{figure}[htbp]
  \centering
  \begin{subfigure}[h]{0.4\textwidth}
    \centering
    \includegraphics[scale = 0.26]{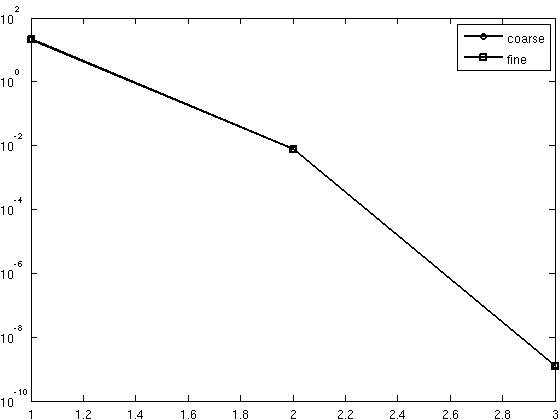}
    \caption{Newton's method}\label{amgmr_newton}
  \end{subfigure}
  \hfill
  \begin{subfigure}[h]{0.4\textwidth}
    \centering
    \includegraphics[scale = 0.25]{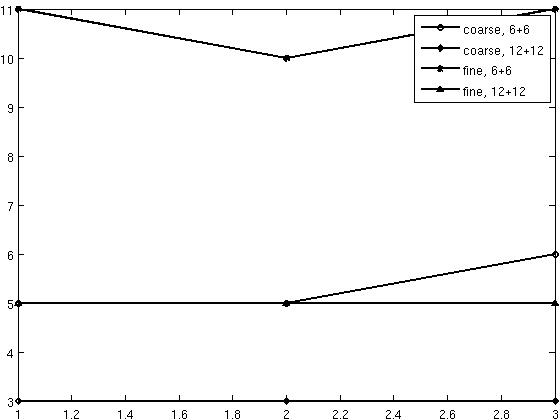}
    \caption{AMG}\label{amgmrr}
  \end{subfigure}
  \caption{Convergence history of Newton's method (left) and 
    iteration numbers of AMG method (right) for the structure sub-problem 
    modeled by the modified hyperelastic Mooney-Rivlin material 
    on coarse and fine meshes, using $6$ and $12$ 
    pre- and post-smoothing steps, respectively.}\label{amgmr}
\end{figure}

For the structure sub-problem modeled by the 
modified hyperelastic Mooney-Rivlin material, 
the convergence of Newton's method and 
the iteration numbers of the AMG method to solve the linearized 
NS equations using different number of {\it Vanka} smoothing steps 
are plotted in Fig. \ref{amgmr}.

As expected, a quadratic convergence rate of Newton's method is observed from the 
numerical results. The number of 
AMG iterations is drastically decreased compared to the preconditioned 
BiCGStab method. In addition, the iteration numbers are less effected by the 
mesh refinement. By doubling the 
pre- and post-smoothing steps (from $6$ to $12$), 
we observe the reduced AMG iteration numbers by almost a factor $2$.

\begin{figure}[htbp]
  \centering
  \begin{subfigure}[h]{0.4\textwidth}
    \centering
    \includegraphics[scale = 0.26]{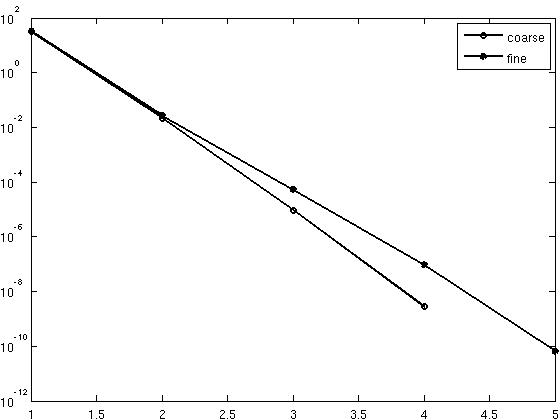}
    \caption{Newton's method}\label{amgart_newton}
  \end{subfigure}
  \hfill
  \begin{subfigure}[h]{0.4\textwidth}
    \centering
    \includegraphics[scale = 0.25]{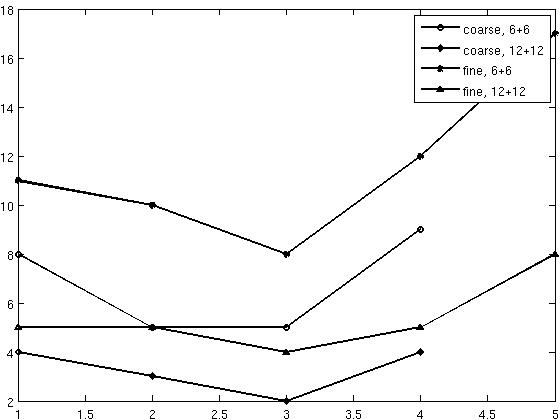}
    \caption{AMG}\label{amgartt}
  \end{subfigure}
  \caption{Convergence history of Newton's method (left) and 
    iteration numbers of AMG method (right) for the structure sub-problem 
    modeled by the two-layer thick-walled artery 
    on coarse and fine meshes, using $6$ and $12$ 
    pre- and post-smoothing steps, respectively.}\label{amgart}
\end{figure}

For the structure sub-problem modeled by the 
two-layer thick-walled artery, the convergence of Newton's method and 
the iteration numbers of the AMG method to solve the linearized 
NS equations using different number of {\it Vanka} smoothing steps 
are plotted in Fig. \ref{amgart}.

As seen from the numerical results, a quadratic convergence rate of Newton's method is 
recovered. The number of 
AMG iterations is drastically decreased compared to the preconditioned 
BiCGStab method. In addition, the iteration numbers are less effected by the 
mesh refinement. By doubling the pre- and post-smoothing steps (from $6$ to $12$), 
we observe the reduced AMG iteration numbers by almost a factor $2$. 

\subsection{Adaptive error control}\label{sec:adpt_err}
For test purpose, we use the adaptive strategy in Section \ref{femnewton} 
to control the error $\varepsilon_2$ of the outer Newton iterations 
and the error $\varepsilon_{1,k}$ of the inner AMG iterations. 
We present the numerical results only for 
the second nonlinear DN FSI iteration 
at the first time step. For other DN iterations and time steps, we 
observe similar results. 

For the fluid sub-problem, we use the AMG method with 
{\it Braess-Sarazin} smoother ($8$ pre- and post-smoothing steps). 
The control of the outer and inner iterations is prescribed 
in Table \ref{tab:adapfluidfinemesh} and \ref{tab:adapfluidfinefinemesh} 
for the coarse and fine meshes, respectively. 
\begin{table}[htbp]
  \begin{tabular}{l|l|l|l}
    \hline
    Newton-It & $\varepsilon_2$& $\varepsilon_{1,k}$ & \#AMG-It \\
    \hline\hline
    $k=1$ & 9.2e+01 & 4.1e-02 & 1 \\
    $k=2$ & 6.2e-03 & 1.8e-03 & 1 \\
    $k=3$ & 6.6e-06 & 1.6e-05 & 2 \\
    $k=4$ & 3.1e-10 & 1.6e-11 & 5 \\
    \hline
\end{tabular}
  \caption{Adaptive error control on the outer (Newton: $\varepsilon_2$) 
    and inner (AMG: $\varepsilon_{1,k}$) iterations 
    for the fluid sub-problem on the coarse mesh.}
  \label{tab:adapfluidfinemesh}
\end{table}
\begin{table}[htbp]
  \begin{tabular}{l|l|l|l}
    \hline
    Newton-It & $\varepsilon_2$& $\varepsilon_{1,k}$ & \#AMG-It \\
    \hline\hline
    $k=1$ &9.4e+01  &7.1e-02 & 1 \\
    $k=2$ &8.6e-03  &1.2e-02 & 1 \\
    $k=3$ &4.4e-05  &4.6e-06 & 3 \\
    $k=4$ &7.9e-10  &1.5e-10 & 6 \\
    \hline
\end{tabular}
  \caption{Adaptive error control on the outer (Newton: $\varepsilon_2$) 
    and inner (AMG: $\varepsilon_{1,k}$) iterations 
    for the fluid sub-problem on the fine mesh.}
  \label{tab:adapfluidfinefinemesh}
\end{table}

It is easy to observe from the numerical results in 
Table \ref{tab:adapfluidfinemesh} and \ref{tab:adapfluidfinefinemesh}, 
the adaptive error control requires fewer inner AMG iterations (\#AMG-It) 
without deterioration of Newton convergence rate of the outer iteration. 

For the structure sub-problem, we use the AMG method with 
{\it Vanka} smoother ($12$ pre- and post-smoothing steps). 
The control of the outer and inner iterations for the 
modified model of Mooney-Rivlin material is shown 
in Table \ref{tab:adapstrucmrfinemesh} 
and \ref{tab:adapstrucmrfinefinemesh} 
for the coarse and fine meshes, respectively.
The control of the outer and inner iterations for the 
two-layer thick-walled arterial model is presented 
in Table \ref{tab:adapstrucgafinemesh} 
and \ref{tab:adapstrucgafinefinemesh} 
for the coarse and fine meshes, respectively. 
\begin{table}[htbp]
  \begin{tabular}{l|l|l|l}
    \hline
    Newton-It  & $\varepsilon_2$& $\varepsilon_{1,k}$ & \#AMG-It \\
    \hline\hline
    $k=1$ &2.4e+01  & 9.2e-07 & 1 \\
    $k=2$ &1.5e-02  & 4.0e-06 & 1 \\
    $k=3$ &1.9e-06  & 1.1e-05 & 2 \\
    $k=4$ &2.6e-11  & 3.3e-13 & 5 \\
    \hline
\end{tabular}
  \caption{Adaptive error control on the outer (Newton: $\varepsilon_2$) 
    and inner (AMG: $\varepsilon_{1,k}$) iterations 
    for the structure sub-problem on the coarse mesh using the 
    modified model of Mooney-Rivlin material.}
  \label{tab:adapstrucmrfinemesh}
\end{table}
\begin{table}[htbp]
  \begin{tabular}{l|l|l|l}
    \hline
    Newton-It & $\varepsilon_2$& $\varepsilon_{1, k}$ & \#AMG-It \\
    \hline\hline
    $k=1$ &2.5e+01  &6.4e-03  & 1 \\
    $k=2$ &1.8e-01  &1.9e-02  & 1 \\ 
    $k=3$ &2.8e-03  &2.8e-02  & 1 \\
    $k=4$ &1.1e-04  &1.4e-06  & 4 \\
    $k=5$ &1.5e-09  &1.8e-09  & 7 \\
    \hline
  \end{tabular}
  \caption{Adaptive error control on the outer (Newton: $\varepsilon_2$) 
    and inner (AMG: $\varepsilon_{1, k}$) iterations 
    for the structure sub-problem on the fine mesh using the 
    modified model of Mooney-Rivlin material.}
  \label{tab:adapstrucmrfinefinemesh}
\end{table}
\begin{table}[htbp]
  \begin{tabular}{l|l|l|l}
    \hline
    Newton-It  & $\varepsilon_2$& $\varepsilon_{1, k}$ & \#AMG-It \\
    \hline\hline
    $k=1$ & 3.2e+01 & 5.0e-06 & 1 \\
    $k=2$ & 2.3e-02 & 1.3e-04 & 1 \\
    $k=3$ & 5.7e-04 & 2.8e-04 & 2 \\
    $k=4$ & 3.0e-06 & 2.0e-07 & 6 \\
    $k=5$ & 1.3e-11 & 2.1e-12 & 9 \\
    \hline
\end{tabular}
  \caption{Adaptive error control on the outer (Newton: $\varepsilon_2$) 
    and inner (AMG: $\varepsilon_{1, k}$) iterations 
    for the structure sub-problem on the coarse mesh using the 
    two-layer thick-walled arterial model.}
  \label{tab:adapstrucgafinemesh}
\end{table}
\begin{table}[htbp]
  \begin{tabular}{l|l|l|l}
    \hline
    Newton-It & $\varepsilon_2$& $\varepsilon_{1, k}$ & \#AMG-It \\
    \hline\hline
    $k=1$ &3.2e+01 &3.2e-05 & 1\\
    $k=2$ &3.6e-02 &3.6e-04 & 1\\
    $k=3$ &3.1e-03 &6.0e-04 & 2\\
    $k=4$ &6.8e-05 &8.3e-06 & 6\\
    $k=5$ &1.5e-09 &4.0e-09 & 9\\
    \hline
  \end{tabular}
  \caption{Adaptive error control on the outer (Newton: $\varepsilon_2$) 
    and inner (AMG: $\varepsilon_{1, k}$) iterations 
    for the structure sub-problem on the fine mesh using the 
    two-layer thick-walled arterial model.}
  \label{tab:adapstrucgafinefinemesh}
\end{table}

From the numerical results in 
Table \ref{tab:adapstrucmrfinemesh}-\ref{tab:adapstrucgafinefinemesh}, 
we observe the satisfying convergence rate of the outer iterations with 
the adaptive error control, that requires 
fewer inner AMG iterations (\#AMG-It). 
The number of Newton iterations is very close to what 
we obtained by using 
the inner AMG method with the fixed error $\varepsilon_1=10^{-8}$.  

\section{Conclusions}\label{sec:con}
We test a partitioned approach for the FSI simulation using two 
hyperelastic models with near-incompressibility constrains. The numerical 
results show the robustness and efficiency of the partitioned 
approach combined with Newton's method to tackle the nonlinear 
fluid and structure sub-problems. Stabilized 
finite element methods, efficient Krylov subspace and AMG 
methods are used to handle the fluid and structure sub-problems, 
that demonstrate the feasibility of our methodology to handle such 
highly nonlinear systems. 
The AMG method shows more robustness 
than the preconditioned Kyrlov subspace method. The adaptive 
error control for the nolinear problems 
requires fewer inner AMG iterations without 
deterioration of the convergence rate of the outer Newton iterations.   
  
\section{Acknowledgement}
We would like to thank Dr. Christoph Augustin and Prof. Gerhard A.
Holzapfel for the discussion on modeling of biological tissues. Special 
thanks go to Prof. Jean-Fr\'{e}d\'{e}ric Gerbeau and Prof. Johan Hoffman 
for the discussion on FSI modeling.
\bibliography{FSI_NonLinea}
\bibliographystyle{siam}

\end{document}